\documentclass[a4paper,fleqn]{cas-sc}

\usepackage[numbers,sort&compress]{natbib}
\usepackage[T1]{fontenc}
\usepackage[utf8]{inputenc}
\renewcommand{\ttdefault}{lmtt}

\usepackage{amsmath,amssymb,amsthm}
\usepackage{booktabs,longtable}
\usepackage{graphicx}
\usepackage{tikz}
\usepackage{xcolor}
\usepackage{url}
\usepackage{algorithm}
\usepackage[noend]{algpseudocode}
\algrenewcommand\algorithmicrequire{\textbf{Input:}}
\algrenewcommand\algorithmicensure{\textbf{Output:}}
\AtBeginDocument{\hypersetup{hidelinks}}

\newtheorem{theorem}{Theorem}
\newtheorem{proposition}[theorem]{Proposition}
\newtheorem{lemma}[theorem]{Lemma}
\newtheorem{corollary}[theorem]{Corollary}
\theoremstyle{definition}
\newtheorem{assumption}{Assumption}
\newtheorem{definition}{Definition}
\newtheorem{example}{Example}

\makeatletter
\def\LT@makecaption#1#2#3{%
  \LT@mcol\LT@cols l{\hbox to\z@{\parbox[t]\linewidth{\raggedright\sffamily\small
    \textbf{\color{scolor}#2}\par#3\par\vskip4pt}\hss}}}
\makeatother

\begin{document}
\let\WriteBookmarks\relax
\shorttitle{Recovering the Reduct Family and the Core from Random Samples}
\shortauthors{Tran}

\title[mode=title]{Recovering the Reduct Family and the Core from Random Samples: Minimax Rates and an Unbounded Separation}

\author[1]{Thanh Dai Tran}[orcid=0009-0003-4527-9754]
\cormark[1]
\credit{Conceptualization, Methodology, Software, Validation, Formal analysis,
  Investigation, Data curation, Writing -- original draft,
  Writing -- review \& editing, Visualization}
\ead{ttdai@uneti.edu.vn}
\affiliation[1]{organization={University of Economics - Technology for Industries},
            addressline={456 Minh Khai, Hai Ba Trung District},
            city={Hanoi},
            postcode={100000},
            country={Vietnam}}
\cortext[1]{Corresponding author}

\begin{highlights}
\item Minimax rates for recovering the reduct family and core on structural-zero classes
\item Two cover radii, computable before sampling, set both rates up to attribute constants
\item Witness-graph vertex covers give exact plug-in error exponents for one population
\item Core recovery is never harder and can need arbitrarily fewer samples than the family
\item Under deterministic labels, an attribute in the empirical core is never spurious
\end{highlights}

\begin{abstract}
Reducts and the core are usually computed from a decision table as if the table were the population. When the table is a random sample, an exactly computed reduct family can differ from the population family, because rare configurations may carry some of its defining distinctions. This paper asks how many observations are needed to recover the population reduct family and its intersection, the core, under finite support and deterministic labels. For model classes given by structural zeros and probability floors, we prove that two weighted cover radii determine the minimax rates of both targets, up to constants depending only on the number of attributes. The radii can be computed before data are collected. Rare conflicts at Hamming distance one limit the core; rare conflicts at any distance limit the family. For a single population, the masses of minimum vertex covers of witness graphs give the exact error exponents of the plug-in estimators. Core recovery is never harder than family recovery, and the gap can be unbounded. In a class of three populations with a rare configuration of mass $\varepsilon$, the core needs a number of observations independent of $\varepsilon$, whereas every family estimator needs at least order $\log (1/\delta )/\varepsilon$. An attribute in the empirical core is never a false positive, whereas the empirical family can stay wrong until its rare witnesses are observed. Exact computation and simulation agree with the formulas, but a study of seven finite empirical populations, specified in advance, found no non-degenerate instance, so prevalence in observational data remains open.
\end{abstract}


\maketitle

\section{Introduction}\label{sec:intro}

Attribute reduction is one of the main uses of rough-set theory. It removes condition attributes that are not needed to determine a decision, and its results are often read as statements about the domain: which attributes are indispensable, and which sets of attributes are enough. Rough-set theory expresses these statements through reducts and the core \cite{pawlak1982,pawlak1991}. In a consistent decision table, a reduct is an inclusion-minimal set of attributes that still determines the decision. The core is the intersection of all reducts. Both are usually computed from the table at hand, under the classical model and under its fuzzy, topological and neighbourhood extensions \cite{dai2024ift,giang2024hausdorff,nguyen2026mhmg}, and the result is taken to describe the population from which the table came.

In practice the table is a random sample, and exact computation does not make the answer correct. A distinction between two differently labelled configurations may be witnessed only by a rare configuration. Until that configuration is observed, the reduct family computed from the sample can differ from the population family, although the computation itself is exact. The difficulty is statistical, not computational. It lies in which distinctions determine the target and how likely a sample is to reveal them.

Within rough-set theory, randomness has entered in several ways. Probabilistic rough sets put probability into the approximation operators \cite{yao2008}, with reducts defined to preserve them \cite{ma2020}, and variable-precision and approximate reducts relax the sufficiency requirement \cite{ziarko1993,slezak2002}, and a variable-precision model has also been built on fuzzy rough sets for noisy numerical tables \cite{ha2024vpfrs}. Other work treats the table itself as data. Incremental methods update reducts when objects are added or removed \cite{kryszkiewicz1997,wei2018,yang2019}, and dynamic reducts measure how often a reduct remains a reduct in random subtables of a fixed table \cite{bazan1994,bazan2000}. The STRIM method treats a decision table as an i.i.d.\@ sample from a population labelled by if-then rules. It derives approximate table sizes for detecting a true rule and keeps the attributes that a $\chi^2$ test finds dependent on the decision \cite{kato2013,kato2014,kato2015,fei2017,kato2019}. Randomisation tests and Hoeffding-type bounds assess the reliability of dependencies and rules \cite{duntsch1997,jaworski2008}, and recent stability-based results bound the generalisation error of a classifier built on a reduct or a granularity \cite{gao2024a,gao2024b}. These methods describe a table, a rule, a dependency or a classifier. They do not address the recovery of the complete reduct family or the core of the population.

Feature selection meets the same objects under other names. With deterministic labels, the core is the set of strongly relevant attributes \cite{kohavi1997} and the reducts are the Markov boundaries \cite{tsamardinos2003} (Section~\ref{sec:population-reducts}). There are scalable algorithms for learning one boundary \cite{pena2007}, algorithms for finding several boundaries \cite{statnikov2013}, and finite-sample guarantees for recovering a single boundary \cite{gao2021}. The minimal-optimal and all-relevant objectives ask for one suitable minimal subset or for every relevant attribute \cite{nilsson2007,gopfert2018}. Learning-theoretic results bound the samples needed to learn concepts that depend on few attributes \cite{almuallim1991,mossel2004}. Testing whether a Boolean function depends on at most $k$ variables can need as many samples as finding those variables \cite{beretta2026}. Other work classifies variables by whether they occur in all, some or none of the admissible subsets \cite{kowal2022}, or ranges a variable's importance over all near-optimal models \cite{fisher2019}. Stability indices measure how much a selected subset varies across samples \cite{nogueira2018}. This work learns, tests or scores individual variables, or recovers one boundary. It does not recover the whole family of minimal sufficient sets and compare it with its intersection.

Database theory studies the same structures as keys. The minimal keys of a relation form a Sperner system, and every Sperner system arises as a key family \cite{demetrovics1979}; the size of representing relations and Armstrong relations have also been studied \cite{tichler2002,beeri1984}. Enumerating keys and deciding key membership raise computational questions, some NP-complete and some still open \cite{lucchesi1978,beeri1979,eiter1995}. Closest in spirit, Kivinen and Mannila sample tuples from a fixed relation to detect dependencies or keys whose error exceeds a threshold \cite{kivinen1994,kivinen1995}. Practical discovery algorithms compare sampled record pairs before validating candidates on the full relation \cite{papenbrock2016}. These guarantees concern approximate validity in a given relation or computational cost, not the exact structure of a population. We are not aware of any result on the number of observations needed to recover the complete reduct family and the core of a population, or of any comparison of the two targets.

Our approach starts from a simple observation. Under deterministic labels, a set of attributes fails to determine the decision only when two configurations with different labels agree on that set. Each reduct and each core attribute is therefore fixed by such conflicting pairs, which we call witnesses, and a sample recovers a target once it contains the witnesses the target needs. How many observations this takes depends on the probability mass of the rarest group of these witnesses, not on how many reducts there are. We develop this idea for finite populations with deterministic labels and make three contributions.

1) We prove minimax rates for both targets on structural-zero classes. Such a class is given by a set of admissible configurations, a subset known to occur, and probability floors for the configurations that occur. Two weighted cover radii determine the rates up to constants that depend only on the number of attributes (Theorem~\ref{thm:structural-zero}). They can be computed before any data are collected, in time polynomial in the numbers of candidate configurations and attributes. For the core, only conflicting pairs that differ in a single attribute matter; for the family, every conflicting pair matters.

2) We give the exact error exponents of the plug-in estimators for a single population. Minimum vertex covers of witness graphs define deletion radii, whose masses fix these exponents (Proposition~\ref{prop:plugin-exponents}) and have a closed form for any number of decision values. The radii also give bounds on both sample complexities (Theorem~\ref{thm:witness-cover}).

3) We show that the two targets can need very different amounts of data. Recovering the core is never harder than recovering the family (Eq.~\eqref{eq:coarsening}), and the gap has no upper bound. We construct a class of three populations indexed by the mass $\varepsilon$ of a rare configuration. On this class, a number of observations that does not depend on $\varepsilon$ suffices for the core, while every family estimator needs at least order $\log (1/\delta )/\varepsilon$ observations (Theorem~\ref{thm:separation}). The two targets also differ in how far a sample can be trusted. A core attribute found in a sample is always a true core attribute (Lemma~\ref{lem:empirical}), while a family found in a sample can remain wrong until its rare witnesses appear.

The separation needs incomplete support. When every combination of attribute values occurs, the reduct is unique, equals the core, and the two problems coincide (Proposition~\ref{prop:full-support}). The rates of the first contribution hold for the worst member of a class, not for a particular population. The paper is theoretical. We check the formulas by exact computation and simulation. We also examined seven finite empirical populations in a study specified in advance, and none showed a non-degenerate separation; how often the separation arises in observational data is therefore open.

The remaining sections build the three results in the order in which they are proved. Section~\ref{sec:prelim} recalls the rough-set notions and the statistical tools used in the proofs. Section~\ref{sec:population} defines the population model, the two targets and the risks that measure how well a sample recovers them. Section~\ref{sec:separation} proves the separation on a class of three populations (the third contribution), because this class shows the mechanism in its simplest form. Section~\ref{sec:witness-cover} turns the mechanism into witness covers for a single population (the second contribution), and Section~\ref{sec:structural-zero} uses these covers to derive the rates on structural-zero classes (the first contribution), which are the most general results. Section~\ref{sec:experiments} checks the formulas by exact computation and simulation and looks for the separation in seven finite empirical populations. Section~\ref{sec:discussion} explains what the results do and do not say, and Section~\ref{sec:conclusion} concludes. Appendix~\ref{app:reproducibility} documents the reproducibility package and Appendix~\ref{app:simulation} the Monte Carlo validation.

\section{Preliminaries}\label{sec:prelim}

This section recalls the rough-set notions used throughout the paper and the statistical tools of the proofs. Section~\ref{sec:population} restates the rough-set notions for a population.

\subsection{Decision tables, reducts and core}\label{sec:decision-systems}

\begin{definition}[decision table]\label{def:table}
A \emph{decision table} is a pair $(U, C\cup\{d\})$, where $U$ is a finite set of objects, $C=\{a_1,\ldots,a_p\}$ is a set of condition attributes and $d$ is a decision attribute, all with finite value ranges. The table is \emph{consistent} if any two objects with the same values on $C$ have the same decision.
\end{definition}

\begin{definition}[reducts and core]\label{def:table-reduct}
Let $(U, C\cup\{d\})$ be a decision table. A subset $B\subseteq C$ is \emph{sufficient} if any two objects that agree on $B$ have the same decision. A \emph{reduct} is an inclusion-minimal sufficient subset, and the \emph{core} is the intersection of all reducts \cite{pawlak1982,pawlak1991}.
\end{definition}

If the table is consistent, $C$ itself is sufficient, so at least one reduct exists. Throughout the paper the labels are deterministic given the full condition vector. Section~\ref{sec:population-model} states this realisability assumption for a population, and Section~\ref{sec:empirical-populations} examines it on the empirical populations. The reduct family is an antichain under set inclusion and may contain many members. It records all minimal sufficient attribute sets, whereas the core records the attributes shared by every such set. The core may be empty even when every reduct is non-empty, so selecting one reduct answers a different question from determining the complete family or its intersection.

\subsection{Discernibility}\label{sec:discernibility}

\begin{definition}[discernibility sets and clutter]\label{def:table-discernibility}
For two objects $u, v$ with different decisions, the \emph{discernibility set} $D(u,v)$ is the set of condition attributes on which they differ \cite{skowron1992}. The family of all such sets is the \emph{discernibility hypergraph} of the table, and its inclusion-minimal edges form its \emph{clutter}.
\end{definition}

A subset of attributes is sufficient if and only if it meets every discernibility set. Hence the reducts are the inclusion-minimal hitting sets, or minimal transversals, of the discernibility hypergraph \cite{vinterbo2000}, and the clutter has the same minimal transversals. Blocker duality recovers the clutter from the reduct family \cite{edmonds1970}. An attribute belongs to the core exactly when its singleton is a discernibility set. Lemma~\ref{lem:discernibility} proves these statements for a population.

Equivalently, the discernibility sets define the clauses of a monotone Boolean function whose prime implicants correspond to the reducts. The reduct family can have exponentially many members, so explicitly enumerating it may require exponential output size. Core membership, in contrast, can be determined by inspecting singleton discernibility sets. Both are computational properties of an observed table, distinct from the statistical question addressed below: how many observations are needed to reveal the witness pairs that determine the population structures.

The following table, with five objects and three binary attributes, illustrates the definitions:
\[
\begin{array}{c|ccc|c}
 & a_1 & a_2 & a_3 & d\\ \hline
u_1 & 0 & 0 & 0 & 0\\
u_2 & 0 & 0 & 1 & 0\\
u_3 & 0 & 1 & 0 & 0\\
u_4 & 1 & 0 & 0 & 1\\
u_5 & 1 & 1 & 1 & 0
\end{array}
\]
Only $u_4$ has decision $1$, so the discernibility sets are $D(u_1,u_4)=\{a_1\}$, $D(u_2,u_4)=\{a_1,a_3\}$, $D(u_3,u_4)=\{a_1,a_2\}$ and $D(u_4,u_5)=\{a_2,a_3\}$. The clutter keeps $\{a_1\}$ and $\{a_2,a_3\}$, because the other two sets contain $\{a_1\}$. Its minimal transversals are the two reducts $\{a_1,a_2\}$ and $\{a_1,a_3\}$, and the core is $\{a_1\}$. The populations of Fig.~\ref{fig:witness} have the same clutter; there, the masses of the configurations decide how hard each target is to recover from a sample.

\subsection{Statistical tools}\label{sec:tools}

The analysis itself draws on established statistical machinery: minimax risk, total variation and two-point testing \cite{lecam1973,lecam1986,tsybakov2009}, tools also used to obtain information-theoretic limits for sparsity recovery \cite{wainwright2009}. Upper bounds follow from explicit estimators. For lower bounds, we construct populations with different targets whose sample laws overlap substantially. Witness graphs connect these probability arguments to the discernibility structure. The simulations in Section~\ref{sec:experiments} report Wilson intervals for Monte Carlo proportions \cite{wilson1927}.

\section{Population rough-set identification}\label{sec:population}

This section restates the rough-set notions of Section~\ref{sec:prelim} for a population, defines the plug-in estimators, and defines the risks and sample complexities that the rest of the paper compares.

\subsection{Population decision model}\label{sec:population-model}

\begin{definition}[population]\label{def:population}
Let $C = \{a_{1}, \dots , a_p\}$ be a finite set of condition attributes, each with a finite value range, and let $Y$ take values in a finite decision space $\mathcal{Y}$. A \emph{population} is a joint law $P$ of $(X_C, Y)$, where $X_C = (X_{a_{1}}, \dots , X_{a_p})$.
\end{definition}

We impose throughout the two assumptions that define the scope of this paper.

\begin{assumption}[finite support]\label{ass:finite-support}
$|\operatorname{supp} (P)| < \infty$, where $\operatorname{supp} (P)$ is the set of labelled atoms $(x, y)$ carrying positive mass.
\end{assumption}

\begin{assumption}[realisability]\label{ass:realisability}
There is a deterministic $f$ with $Y = f(X_C)$ $P$-almost surely; equivalently $\mathrm{H}_P(Y \mid X_C) = 0$.
\end{assumption}

Since the value ranges are finite, Assumption~\ref{ass:finite-support} holds automatically; it is stated separately because it is the property the arguments use. Realisability is a substantive assumption. It says that the decision is a function of the full condition vector on the support, so the population contains no two identical condition vectors with different labels. Section~\ref{sec:empirical-populations} reports how often it fails on real decision tables.

Under Assumptions~\ref{ass:finite-support} and~\ref{ass:realisability} a population has two layers. Its structure, namely which attribute subsets are sufficient and which of them are minimal, depends only on the \emph{labelled support}: which atoms occur and with which labels. The difficulty of identifying that structure from random observations depends on the labelled support and on its probability masses. Two populations with the same labelled support but different masses therefore have the same rough-set structures, yet they can be arbitrarily far apart in sampling difficulty. Sections~\ref{sec:separation}--\ref{sec:structural-zero} quantify this gap.

\subsection{Population reduct family and core}\label{sec:population-reducts}

\begin{definition}[population reducts and core]\label{def:reduct}
Let $P$ be a population. A subset $B \subseteq C$ is \emph{sufficient} if $Y$ is a deterministic function of $X_B$ almost surely, equivalently $\mathrm{H}_P(Y \mid X_B) = 0$; write $\mathcal{S}_P$ for the family of sufficient subsets. A \emph{population reduct} is an inclusion-minimal element of $\mathcal{S}_P$. The \emph{population reduct family} and the \emph{population core} are
\begin{equation}\label{eq:reduct-core}
\mathcal{R}_P = \{ R \subseteq C: R \ \text{is a population reduct}\  \},\quad \mathrm{CORE}_P = \bigcap _{R \in \mathcal{R}_P} R.
\end{equation}
\end{definition}

By Assumption~\ref{ass:realisability}, $C \in \mathcal{S}_P$, and $\mathcal{S}_P$ is upward closed, so $\mathcal{R}_P$ is non-empty. These are the classical objects of rough-set theory at the population level, and we use one standard representation of them throughout.

\begin{definition}[discernibility sets and clutter]\label{def:discernibility}
For oppositely labelled support points $x, x'$ of $P$, the \emph{discernibility set} is $D(x, x') = \{a \in C: x_a \neq x'_a\}$. Let $\mathcal{H}_P$ be the family of all such discernibility sets, and let $\mathcal{C}_P = \min _{\subseteq} \mathcal{H}_P$ be its \emph{clutter} of inclusion-minimal edges.
\end{definition}

\begin{lemma}[discernibility representation]\label{lem:discernibility}
Under Assumptions~\ref{ass:finite-support} and~\ref{ass:realisability}:
(i) $B$ is sufficient if and only if $B \cap e \neq \emptyset$ for every $e \in \mathcal{H}_P$;
(ii) $\mathcal{R}_P = \mathrm{Tr}_{\mathrm{min}}(\mathcal{H}_P) = \mathrm{Tr}_{\mathrm{min}}(\mathcal{C}_P)$, the minimal transversals;
(iii) $\mathrm{CORE}_P = \{a \in C: \{a\} \in \mathcal{C}_P\}$.
\end{lemma}

\begin{proof}

(i) $B$ fails to be sufficient exactly when two oppositely labelled support points agree on $B$, and $x_B = x'_B$ is equivalent to $B \cap D(x, x') = \emptyset$. (ii) is immediate from (i) and the fact that every edge contains a minimal edge. (iii) If $\{a\} \in \mathcal{H}_P$ then every transversal contains $a$, so $a$ lies in every reduct. Conversely if $\{a\} \notin \mathcal{H}_P$ then every edge meets $C \setminus \{a\}$, so $C \setminus \{a\}$ is sufficient and therefore contains a reduct avoiding $a$.
\end{proof}

Lemma~\ref{lem:discernibility} is classical: (i)--(ii) are the discernibility machinery of Skowron and Rauszer~\cite{skowron1992}, restated at the population level, and (iii) is the standard characterisation of the core by singleton edges. The proof fixes notation for the empirical analogue in Section~\ref{sec:empirical}; blocker duality \cite{edmonds1970} will identify the family with its clutter.

For later interpretation, writing $Z_a = X_{C \setminus \{a\}}$ and $\tau _a = \inf _g P\{Y \neq g(Z_a)\}$, one has $a \in \mathrm{CORE}_P$ if and only if $\tau _a > 0$; we call this the leave-one-out characterisation. Equivalently, under Assumption~\ref{ass:realisability} the population core is exactly the set of \emph{strongly relevant} attributes in the sense of Kohavi and John~\cite{kohavi1997}. The same assumption identifies the reducts with Markov boundaries \cite{tsamardinos2003}. A set $B\subseteq C$ is a Markov blanket when $Y$ is conditionally independent of $X_{C\setminus B}$ given $X_B$, and an inclusion-minimal blanket is a Markov boundary. Under Assumption~\ref{ass:realisability}, the blanket condition forces $P(Y\mid X_B)=P(Y\mid X_C)$ almost surely, and the latter is a point mass, so $B$ is sufficient. Conversely, if $B$ is sufficient, then $Y$ is determined by $X_B$, which implies the blanket condition. The reducts are therefore the Markov boundaries, and the core is the intersection of all boundaries. These equivalences connect the recovery problem to feature relevance.

\subsection{Empirical structures}\label{sec:empirical}

\begin{definition}[empirical reducts and core]\label{def:empirical}
Let $S_n = ((X_{1}, Y_{1}), \dots , (X_n, Y_n))$ be i.i.d.\ from $P$. A subset $B \subseteq C$ is \emph{empirically sufficient} if no two observations agree on $B$ and disagree on the label. The \emph{empirical reducts} are the inclusion-minimal empirically sufficient sets; $\widehat{\mathcal{R}}_n$ is the complete family of them, and $\widehat{\mathrm{CORE}}_n = \bigcap _{R \in \widehat{\mathcal{R}}_n} R$. We write $\widehat{\mathcal{H}}_n$ for the family of discernibility sets of oppositely labelled pairs both of whose atoms appear in the sample.
\end{definition}

The estimators $\widehat{\mathcal{R}}_n$ and $\widehat{\mathrm{CORE}}_n$ are the \emph{plug-in estimators} of the two targets. Duplicated observations are irrelevant to Definition~\ref{def:empirical}: $\widehat{\mathcal{R}}_n$ and $\widehat{\mathrm{CORE}}_n$ depend on the sample only through the \emph{set of distinct labelled atoms observed}, and the proof of Lemma~\ref{lem:discernibility} applies verbatim to that set.

\begin{lemma}[empirical representation and one-sidedness]\label{lem:empirical}
Under Assumptions~\ref{ass:finite-support} and~\ref{ass:realisability}, for every sample:
(i) $\widehat{\mathcal{R}}_n = \mathrm{Tr}_{\mathrm{min}}(\widehat{\mathcal{H}}_n)$ and $\widehat{\mathrm{CORE}}_n = \{a: \{a\} \in \widehat{\mathcal{H}}_n\}$;
(ii) $\widehat{\mathcal{H}}_n \subseteq \mathcal{H}_P$, and hence $\widehat{\mathrm{CORE}}_n \subseteq \mathrm{CORE}_P$;
(iii) the plug-in estimators are correct exactly when the sample realises the following edges:
\begin{gather}
\widehat{\mathrm{CORE}}_n = \mathrm{CORE}_P \iff \ \text{every singleton edge of}\  \mathcal{C}_P \ \text{is realised by an observed pair}\ ; \label{eq:core-recovery}\\
\widehat{\mathcal{R}}_n = \mathcal{R}_P \iff \ \text{every edge of}\  \mathcal{C}_P \ \text{is realised by an observed pair}\ . \label{eq:family-recovery}
\end{gather}
\end{lemma}

\begin{proof}
(i) is Lemma~\ref{lem:discernibility} applied to the observed labelled support. (ii) Every observed pair is a population pair, so $\widehat{\mathcal{H}}_n \subseteq \mathcal{H}_P$. If $\{a\} \in \widehat{\mathcal{H}}_n$, then $\{a\} \in \mathcal{H}_P$, so $\{a\} \in \mathcal{C}_P$, because a singleton is minimal wherever it appears, and $a \in \mathrm{CORE}_P$ by Lemma~\ref{lem:discernibility}(iii).

(iii) By (i), (ii) and Lemma~\ref{lem:discernibility}(iii), $\widehat{\mathrm{CORE}}_n = \{a: \{a\} \in \widehat{\mathcal{H}}_n \cap \mathcal{C}_P\}$, which equals $\mathrm{CORE}_P = \{a: \{a\} \in \mathcal{C}_P\}$ exactly when every singleton edge of $\mathcal{C}_P$ lies in $\widehat{\mathcal{H}}_n$. This is Eq.~\eqref{eq:core-recovery}. For Eq.~\eqref{eq:family-recovery}, write $\widehat{\mathcal{C}}_n = \min_{\subseteq} \widehat{\mathcal{H}}_n$, so that $\widehat{\mathcal{R}}_n = \mathrm{Tr}_{\mathrm{min}}(\widehat{\mathcal{C}}_n)$. If $\mathcal{C}_P \subseteq \widehat{\mathcal{H}}_n$, then
\[
\mathcal{C}_P \subseteq \widehat{\mathcal{H}}_n \subseteq \mathcal{H}_P ,
\]
and every element of $\mathcal{H}_P$ contains an element of $\mathcal{C}_P$; hence $\widehat{\mathcal{C}}_n = \mathcal{C}_P$ and $\widehat{\mathcal{R}}_n = \mathcal{R}_P$. If some $e \in \mathcal{C}_P$ satisfies $e \notin \widehat{\mathcal{H}}_n$, then $e \notin \widehat{\mathcal{C}}_n$, so $\widehat{\mathcal{C}}_n \neq \mathcal{C}_P$; by blocker duality \cite{edmonds1970} two different clutters have different families of minimal transversals, so $\widehat{\mathcal{R}}_n \neq \mathcal{R}_P$. If no oppositely labelled pair has been observed, then $\widehat{\mathcal{H}}_n = \emptyset$, $\widehat{\mathcal{R}}_n = \{\emptyset\}$ and $\widehat{\mathrm{CORE}}_n = \emptyset$, and both criteria still apply.
\end{proof}

By Lemma~\ref{lem:empirical}(ii), under Assumption~\ref{ass:realisability} sampling can never produce a false core attribute. Part (iii) turns both recovery events into statements about which atoms have been observed; Section~\ref{sec:experiments} uses it throughout.

\subsection{Exact-recovery risks}\label{sec:risks}

This subsection defines the risks compared throughout the paper.

\begin{definition}[recovery risks]\label{def:risks}
Let $\Psi$ be a target functional, here $\Psi = \mathrm{CORE}$ or $\Psi = \mathcal{R}$, and let $\mathcal{P}$ be a class of populations.
\begin{enumerate}
\item[(a)] The \emph{minimax risk} is
\begin{equation}\label{eq:minimax-risk}
\mathfrak{R}_n(\Psi ; \mathcal{P}) = \inf _{\hat{\Psi}_n} \sup _{P \in \mathcal{P}} P^{n}\{ \hat{\Psi}_n \neq \Psi (P) \},
\end{equation}
the infimum running over all estimators measurable with respect to the sample and an independent randomisation, with probability also taken over that randomisation.
\item[(b)] The \emph{plug-in risk} at $P$ is $P^{n}\{ \widehat{\mathrm{CORE}}_n \neq \mathrm{CORE}_P \}$ for the core and $P^{n}\{ \widehat{\mathcal{R}}_n \neq \mathcal{R}_P \}$ for the family, with the plug-in estimators of Definition~\ref{def:empirical}.
\item[(c)] The \emph{two-point testing risk} of $P, Q \in \mathcal{P}$ is the Bayes risk of deciding between $P$ and $Q$ from the sample under a uniform prior.
\end{enumerate}
\end{definition}

The supremum of the plug-in risk over $\mathcal{P}$ upper-bounds the minimax risk and is generally larger. The threshold ratio reported in Section~\ref{sec:empirical-populations} is a descriptive statistic, not one of these risks.

All minimax comparisons below use the same class and confidence level. An estimator may depend on that class, but not on which member generated the observations.

\subsection{Identification sample complexity}\label{sec:sample-complexity}

\begin{definition}[identification sample complexity]\label{def:sample-complexity}
For $0 < \delta < 1$ and a class $\mathcal{P}$, the core and family sample complexities are
\begin{gather}
n_{\mathrm{core}}(\delta ; \mathcal{P}) = \min \{ n \ge 1: \mathfrak{R}_n(\mathrm{CORE}; \mathcal{P}) \le \delta \}, \label{eq:n-core}\\
n_{\mathrm{family}}(\delta ; \mathcal{P}) = \min \{ n \ge 1: \mathfrak{R}_n(\mathcal{R}; \mathcal{P}) \le \delta \}, \label{eq:n-family}
\end{gather}
with the minimum equal to $+\infty$ if the set is empty.
\end{definition}

Both targets use the same labelled sample, class knowledge and independent randomisation. At each fixed finite-support population, the plug-in estimates are eventually correct almost surely: all support atoms are eventually observed. This pointwise consistency does not give a sample bound uniform over populations with arbitrarily small informative masses.

One degeneracy has to be excluded. If the class $\mathcal{P}$ is chosen so small that $\mathrm{CORE}_P$ is constant on it, then the constant estimator that ignores the data identifies the core with zero error and $n_{\mathrm{core}}(\delta ; \mathcal{P}) = 1$ for every $\delta$. The comparison then becomes vacuous, and any ratio it produces is an artefact of the class rather than a statement about identification. For this reason the class used in Section~\ref{sec:separation} contains three populations, on which \emph{both} functionals are non-constant. Non-constancy is necessary but not sufficient: a class can make both functionals non-constant while its core-changing member is so far away that the core is identified from a single observation. Section~\ref{sec:witness-cover} therefore compares the two functionals on classes whose members are the population's own nearest family-changing and core-changing deletions, and Section~\ref{sec:structural-zero} on classes fixed in advance by structural zeros and probability floors.

\subsection{The central comparison}\label{sec:comparison}

The paper compares the two sample complexities on a common class. Fix a class $\mathcal{P}$ on which both functionals are non-constant, and fix $\delta$. The question is how large the ratio

\begin{equation}\label{eq:ratio}
n_{\mathrm{family}}(\delta ; \mathcal{P}) / n_{\mathrm{core}}(\delta ; \mathcal{P})
\end{equation}

can become. If it were bounded by a universal constant, recovering the core and recovering the family would be statistically equivalent up to a constant factor. It would then be enough to study either one. Section~\ref{sec:separation} shows that no such constant exists.

The ratio is never below one. Let $\tilde{\mathcal{R}}_n$ be any family estimator and put $\tilde{C}_n = \bigcap \tilde{\mathcal{R}}_n$. Since $\mathrm{CORE}_P = \bigcap \mathcal{R}_P$, we have $\{\tilde{C}_n \neq \mathrm{CORE}_P\} \subseteq \{\tilde{\mathcal{R}}_n \neq \mathcal{R}_P\}$ for every $P \in \mathcal{P}$. Taking the supremum over $\mathcal{P}$ and then the infimum over $\tilde{\mathcal{R}}_n$ gives
\begin{equation}\label{eq:coarsening}
\mathfrak{R}_n(\mathrm{CORE}; \mathcal{P}) \le \mathfrak{R}_n(\mathcal{R}; \mathcal{P}) \quad\text{and}\quad n_{\mathrm{core}}(\delta ; \mathcal{P}) \le n_{\mathrm{family}}(\delta ; \mathcal{P})
\end{equation}
on every class. Section~\ref{sec:separation} gives an explicit separation within realisable Pawlak populations, and Section~\ref{sec:witness-cover} a witness-cover condition that produces further separating classes.

\section{Core--family separation on a three-population class}\label{sec:separation}

\subsection{Construction and recovery witnesses}\label{sec:construction}

Fix $p = 3$ with attributes named $a_{1}$, $a_{2}$, $a_{3}$, so that $C = \{a_{1}, a_{2}, a_{3}\}$. Take four atoms of $\{0,1\}^{3}$,

\begin{equation}\label{eq:atoms}
x_{0} = (0,0,0),\quad x_{1} = (0,0,1),\quad x_{2} = (1,0,0),\quad x_{3} = (1,1,1),
\end{equation}

and for $0 < \varepsilon \le 1/4$ the masses

\begin{equation}\label{eq:masses}
P(x_{0}) = 1/4,\quad P(x_{1}) = 1/2 - \varepsilon ,\quad P(x_{2}) = 1/4,\quad P(x_{3}) = \varepsilon .
\end{equation}

Three label assignments give three populations, collected in the class $\mathcal{P}_{\varepsilon} = \{P_{\varepsilon}^{(0)}, P_{\varepsilon}^{(1)}, P_{\varepsilon}^{(2)}\}$. Table~\ref{tab:populations} records the assignments together with the exact structures they induce.

\begin{table}[pos=!htbp]
\centering
\small
\caption{The three populations of $\mathcal{P}_\varepsilon$: masses, label assignments and the structures they induce.}\label{tab:populations}
\begin{tabular}{llllllll}
\toprule
 & $x_{0}$ & $x_{1}$ & $x_{2}$ & $x_{3}$ & clutter $\mathcal{C}_P$ & $\mathcal{R}_P$ & $\mathrm{CORE}_P$ \\
\midrule
mass & $1/4$ & $1/2 - \varepsilon$ & $1/4$ & $\varepsilon$ &  &  &  \\
$P_{\varepsilon}^{(0)}$ & 0 & 0 & 1 & 0 & $\{a_{1}\}, \{a_{2},a_{3}\}$ & $\{\{a_{1},a_{2}\}, \{a_{1},a_{3}\}\}$ & $\{a_{1}\}$ \\
$P_{\varepsilon}^{(1)}$ & 0 & 0 & 1 & 1 & $\{a_{1}\}$ & $\{\{a_{1}\}\}$ & $\{a_{1}\}$ \\
$P_{\varepsilon}^{(2)}$ & 0 & 1 & 0 & 1 & $\{a_{3}\}$ & $\{\{a_{3}\}\}$ & $\{a_{3}\}$ \\
\bottomrule
\end{tabular}
\end{table}

All three populations are realisable, since the labels are functions of the atom. The core functional is non-constant on $\mathcal{P}_{\varepsilon}$ (it is $\{a_{1}\}$ twice and $\{a_{3}\}$ once) and so is the reduct family, so the degeneracy of Section~\ref{sec:sample-complexity} does not arise.

\textbf{Discernibility calculation.}

The structures follow from Lemma~\ref{lem:discernibility} by listing conflicting pairs, and the listing is short enough to verify by hand.

Under $P_{\varepsilon}^{(0)}$ the only atom labelled 1 is $x_{2}$, so the conflicting pairs are $(x_{0}, x_{2}) \to \{a_{1}\}$, $(x_{1}, x_{2}) \to \{a_{1}, a_{3}\}$ and $(x_{2}, x_{3}) \to \{a_{2}, a_{3}\}$. Every transversal must contain $a_{1}$ because of the singleton edge; $a_{1}$ already meets $\{a_{1},a_{3}\}$; and $\{a_{2},a_{3}\}$ forces one of $a_{2}$ or $a_{3}$. Hence $\mathcal{R}_{P^{(0)}} = \{\{a_{1},a_{2}\}, \{a_{1},a_{3}\}\}$ and $\mathrm{CORE}_{P^{(0)}} = \{a_{1}\}$.

Under $P_{\varepsilon}^{(1)}$ the atoms labelled 1 are $x_{2}$ and $x_{3}$, giving $(x_{0},x_{2}) \to \{a_{1}\}$, $(x_{0},x_{3}) \to \{a_{1},a_{2},a_{3}\}$, $(x_{1},x_{2}) \to \{a_{1},a_{3}\}$ and $(x_{1},x_{3}) \to \{a_{1},a_{2}\}$. Every edge contains $a_{1}$, so $\{a_{1}\}$ is the unique minimal transversal: $\mathcal{R}_{P^{(1)}} = \{\{a_{1}\}\}$ and $\mathrm{CORE}_{P^{(1)}} = \{a_{1}\}$.

Under $P_{\varepsilon}^{(2)}$ the label is the value of $a_{3}$, and the pairs are $(x_{0},x_{1}) \to \{a_{3}\}$, $(x_{0},x_{3}) \to \{a_{1},a_{2},a_{3}\}$, $(x_{2},x_{1}) \to \{a_{1},a_{3}\}$ and $(x_{2},x_{3}) \to \{a_{2},a_{3}\}$. The singleton $\{a_{3}\}$ meets every edge, so $\mathcal{R}_{P^{(2)}} = \{\{a_{3}\}\}$ and $\mathrm{CORE}_{P^{(2)}} = \{a_{3}\}$.

The comparison of the first two lines is the mechanism of the whole paper. Passing from $P_{\varepsilon}^{(1)}$ to $P_{\varepsilon}^{(0)}$, which flips the label of the single rare atom $x_{3}$, deletes the edges $\{a_{1},a_{2}\}$ and $\{a_{1},a_{2},a_{3}\}$ and creates the edge $\{a_{2},a_{3}\}$. The singleton edge $\{a_{1}\}$ is untouched, so by Lemma~\ref{lem:discernibility}(iii) the core does not move; but the clutter changes, so the family of minimal transversals does. Two populations therefore share a core and differ in their reduct family, and the entire difference between them is carried by an atom of mass $\varepsilon$.

\textbf{Recovery witnesses.}

Lemma~\ref{lem:empirical} turns both recovery events into observation events, and on $\mathcal{P}_{\varepsilon}$ the resulting witness sets can be written down completely (Table~\ref{tab:witnesses}).

\begin{table}[pos=!htbp]
\centering
\small
\caption{Recovery witnesses for the populations of Table~\ref{tab:populations}: the atoms whose observation is equivalent to exact recovery of each target.}\label{tab:witnesses}
\begin{tabular}{lll}
\toprule
population & core recovered iff observed & complete family recovered iff observed \\
\midrule
$P_{\varepsilon}^{(0)}$ & $\{x_{0}, x_{2}\}$ & $\{x_{0}, x_{2}, x_{3}\}$ \\
$P_{\varepsilon}^{(1)}$ & $\{x_{0}, x_{2}\}$ & $\{x_{0}, x_{2}\}$ \\
$P_{\varepsilon}^{(2)}$ & $\{x_{0}, x_{1}\}$ & $\{x_{0}, x_{1}\}$ \\
\bottomrule
\end{tabular}
\end{table}

Exact recovery of either structure is equivalent to having observed a specific finite set of atoms, and under $P_{\varepsilon}^{(0)}$ the two requirements differ by the single atom $x_{3}$ of mass $\varepsilon$. The structural difference is one atom; the statistical difference is a waiting time of order $1/\varepsilon$. Under $P_{\varepsilon}^{(0)}$ and $P_{\varepsilon}^{(1)}$, both core witness atoms have mass $1/4$; under $P_{\varepsilon}^{(2)}$, both have mass at least $1/4$. This gives the uniform core bound below.

Call a set of atoms a \emph{certificate} for a target if observing all of its atoms guarantees exact plug-in recovery of that target; the witness sets of Table~\ref{tab:witnesses} are the minimal certificates on $\mathcal{P}_{\varepsilon}$. In general a minimal certificate chooses one witness pair for each relevant clutter edge, and Section~\ref{sec:witness-cover} turns the resulting waiting times into error exponents. On $\mathcal{P}_{\varepsilon}$ every clutter edge has a single witness pair (Fig.~\ref{fig:witness}a), which is why the witness sets above are unique.

\subsection{Core identification is uniform in $\varepsilon$}\label{sec:core-uniform}

By Lemma~\ref{lem:empirical}, $\widehat{\mathrm{CORE}}_n$ equals the set of attributes appearing as observed singleton edges. Under $P_{\varepsilon}^{(0)}$ and $P_{\varepsilon}^{(1)}$ the only pair with a singleton discernibility set is $(x_{0}, x_{2})$, and under $P_{\varepsilon}^{(2)}$ it is $(x_{0}, x_{1})$. Core recovery is therefore exactly the event that a specific pair of atoms has been observed, and inclusion--exclusion gives closed forms.

\begin{proposition}[exact core recovery]\label{prop:core-exact}
For every $n \ge 1$ and every $0 < \varepsilon \le 1/4$,
\begin{equation}\label{eq:core-exact-01}
P^{(0)}\{ \widehat{\mathrm{CORE}}_n = \mathrm{CORE} \} = P^{(1)}\{ \widehat{\mathrm{CORE}}_n = \mathrm{CORE} \} = 1 - 2(3/4)^{n} + (1/2)^{n},
\end{equation}
a quantity independent of $\varepsilon$, and
\begin{equation}\label{eq:core-exact-2}
e_{\mathrm{core}}^{(2)}(n,\ \varepsilon ):= P^{(2)}\{ \widehat{\mathrm{CORE}}_n \neq \mathrm{CORE} \} = (3/4)^{n} + (1/2 + \varepsilon )^{n} - (1/4 + \varepsilon )^{n}.
\end{equation}
\end{proposition}

\begin{proof}
For an atom $x$ write $N_x$ for the event that $x$ is observed at least once among the $n$ draws. For a set $A$ of atoms, the probability that no atom of $A$ is observed is $(1 - P(A))^{n}$. Under $P^{(0)}$ and $P^{(1)}$, Lemma~\ref{lem:empirical}(iii) and Table~\ref{tab:witnesses} give $\{\widehat{\mathrm{CORE}}_n = \mathrm{CORE}\} = N_{x_{0}} \cap N_{x_{2}}$, and $P(x_{0}) = P(x_{2}) = 1/4$. By inclusion--exclusion,
\begin{align*}
P(N_{x_{0}} \cap N_{x_{2}}) &= 1 - P(N_{x_{0}}^{c}) - P(N_{x_{2}}^{c}) + P(N_{x_{0}}^{c} \cap N_{x_{2}}^{c})\\
&= 1 - 2(1 - 1/4)^{n} + (1 - 1/2)^{n} = 1 - 2(3/4)^{n} + (1/2)^{n},
\end{align*}
which is Eq.~\eqref{eq:core-exact-01}. Under $P^{(2)}$, the same results give $\{\widehat{\mathrm{CORE}}_n = \mathrm{CORE}\} = N_{x_{0}} \cap N_{x_{1}}$, with $P(x_{0}) = 1/4$ and $P(x_{1}) = 1/2 - \varepsilon$. Hence
\begin{align*}
e_{\mathrm{core}}^{(2)}(n, \varepsilon) = P(N_{x_{0}}^{c} \cup N_{x_{1}}^{c}) &= P(N_{x_{0}}^{c}) + P(N_{x_{1}}^{c}) - P(N_{x_{0}}^{c} \cap N_{x_{1}}^{c})\\
&= (3/4)^{n} + (1/2 + \varepsilon)^{n} - (1/4 + \varepsilon)^{n},
\end{align*}
which is Eq.~\eqref{eq:core-exact-2}.
\end{proof}

\begin{corollary}[uniform core bound]\label{cor:core-uniform}
For $0 < \varepsilon \le 1/4$ and all $n$,
\begin{equation}\label{eq:core-uniform}
\sup _{P \in \mathcal{P}_{\varepsilon}} P^{n}\{ \widehat{\mathrm{CORE}}_n \neq \mathrm{CORE}_P \} \le 2(3/4)^{n},
\end{equation}
and hence $n_{\mathrm{core}}(\delta ; \mathcal{P}_{\varepsilon}) \le \lceil \log (2/\delta ) / \log (4/3) \rceil$, an integer independent of $\varepsilon$.
\end{corollary}

\begin{proof}
For $P \in \{P^{(0)}, P^{(1)}\}$, Eq.~\eqref{eq:core-exact-01} gives
\[
P^{n}\{\widehat{\mathrm{CORE}}_n \neq \mathrm{CORE}_P\} = 2(3/4)^{n} - (1/2)^{n} \le 2(3/4)^{n}.
\]
For $P^{(2)}$, the assumption $\varepsilon \le 1/4$ gives $1/2 + \varepsilon \le 3/4$, and $(1/4 + \varepsilon)^{n} \ge 0$; hence, by Eq.~\eqref{eq:core-exact-2},
\[
e_{\mathrm{core}}^{(2)}(n, \varepsilon) \le (3/4)^{n} + (3/4)^{n} = 2(3/4)^{n}.
\]
This proves Eq.~\eqref{eq:core-uniform}. The minimax risk is at most the supremum of the plug-in risk, so $\mathfrak{R}_n(\mathrm{CORE}; \mathcal{P}_{\varepsilon}) \le 2(3/4)^{n} \le \delta$ whenever $n \ge \log (2/\delta) / \log (4/3)$, and the bound on $n_{\mathrm{core}}$ follows from Eq.~\eqref{eq:n-core}.
\end{proof}

Uniformity holds because every required core witness atom has mass at least $1/4$, including $x_{1}$ under $P_{\varepsilon}^{(2)}$ when $\varepsilon \le 1/4$. The core problem is nonetheless non-trivial, which is what makes the comparison in Section~\ref{sec:main-theorem} meaningful.

\begin{proposition}[core lower bound]\label{prop:core-lower}
For $0 < \varepsilon \le 1/4$,
\begin{equation}\label{eq:core-lower}
\mathfrak{R}_n(\mathrm{CORE}; \mathcal{P}_{\varepsilon}) \ge \tfrac{1}{2} (1/4 + \varepsilon )^{n} \ge \tfrac{1}{2} (1/4)^{n},\quad \ \text{hence}\  n_{\mathrm{core}}(\delta ; \mathcal{P}_{\varepsilon}) \ge \log (1/(2\delta )) / \log 4.
\end{equation}
\end{proposition}

\begin{proof}
Consider the subclass $\{P^{(1)}, P^{(2)}\}$, whose cores are $\{a_{1}\}$ and $\{a_{3}\}$ (Table~\ref{tab:populations}). Let $Q_{1}$ and $Q_{2}$ be the laws of one labelled observation under $P^{(1)}$ and $P^{(2)}$, and write $Q_{1} \wedge Q_{2} = \sum_{z} \min (Q_{1}(z), Q_{2}(z))$ for their overlap. By Table~\ref{tab:populations}, the labelled atoms with positive mass under both laws are $(x_{0}, 0)$, of mass $1/4$, and $(x_{3}, 1)$, of mass $\varepsilon$; the atoms $x_{1}$ and $x_{2}$ carry different labels. Hence
\[
Q_{1} \wedge Q_{2} = 1/4 + \varepsilon , \qquad Q_{1}^{\otimes n} \wedge Q_{2}^{\otimes n} \ge (Q_{1} \wedge Q_{2})^{n} = (1/4 + \varepsilon)^{n},
\]
where the inequality uses $\min (ab, cd) \ge \min (a, c) \min (b, d)$ for non-negative $a, b, c, d$. Let $\widehat{C}_n$ be any core estimator, randomised or not, and define the test $\psi = 2$ if $\widehat{C}_n = \{a_{3}\}$ and $\psi = 1$ otherwise. Then $\{\psi \neq k\} \subseteq \{\widehat{C}_n \neq \mathrm{CORE}_{P^{(k)}}\}$ for $k \in \{1, 2\}$, and the two-point bound \cite{lecam1986,tsybakov2009} gives
\[
\max_{k \in \{1, 2\}} P^{(k)}\{\widehat{C}_n \neq \mathrm{CORE}_{P^{(k)}}\} \ge \tfrac{1}{2} \sum_{k = 1}^{2} P^{(k)}\{\psi \neq k\} \ge \tfrac{1}{2}\, Q_{1}^{\otimes n} \wedge Q_{2}^{\otimes n} \ge \tfrac{1}{2} (1/4 + \varepsilon)^{n}.
\]
The supremum over $\mathcal{P}_{\varepsilon}$ is at least the maximum over the subclass, which proves the first inequality of Eq.~\eqref{eq:core-lower}; the second uses $\varepsilon > 0$. If $n < \log (1/(2\delta)) / \log 4$, then $\tfrac{1}{2}(1/4)^{n} > \delta$, so $\mathfrak{R}_n(\mathrm{CORE}; \mathcal{P}_{\varepsilon}) > \delta$; this gives the bound on $n_{\mathrm{core}}$.
\end{proof}

The sample-complexity bounds in Eqs.~\eqref{eq:core-uniform} and~\eqref{eq:core-lower} are independent of $\varepsilon$, so $n_{\mathrm{core}}(\delta ; \mathcal{P}_{\varepsilon}) = \Theta (\log (1/\delta ))$ uniformly in $\varepsilon$. At $\delta = 0.01$ they read $3 \le n_{\mathrm{core}} \le 19$. The upper bound is attained by the class-agnostic plug-in estimator; an estimator that knows the class does much better, since $P^{(2)}$ is revealed by any observation of $(x_{1}, 1)$ or $(x_{2}, 0)$. The value $19$ is therefore $n_{\mathrm{core}}^{\mathrm{UB}}$, an upper bound, not the core sample complexity itself.

\subsection{The complete-family obstruction}\label{sec:family-obstruction}

The family target behaves differently, and the reason is visible in Table~\ref{tab:populations}: $P_{\varepsilon}^{(0)}$ and $P_{\varepsilon}^{(1)}$ are distinguished only by the label attached to an atom of mass $\varepsilon$.

\begin{lemma}[observational indistinguishability]\label{lem:indistinguishability}
Let $E_n$ be the event that no observation has $X = x_{3}$, and write $\mathcal{L}_P(S_n \mid E_n)$ for the conditional law of the labelled sample $S_n$ given $E_n$ under $P$. Then
\begin{equation}\label{eq:indist}
P^{(0)}(E_n) = P^{(1)}(E_n) = (1 - \varepsilon )^{n}, \qquad \mathcal{L}_{P^{(0)}}(S_n \mid E_n) = \mathcal{L}_{P^{(1)}}(S_n \mid E_n).
\end{equation}
\end{lemma}

\begin{proof}
Let $f_{0}$ and $f_{1}$ be the labellings of $P^{(0)}$ and $P^{(1)}$ in Table~\ref{tab:populations}. Both populations have the masses of Eq.~\eqref{eq:masses}, so, by independence of the draws,
\[
P^{(0)}(E_n) = \prod_{i=1}^{n} P(X_i \neq x_{3}) = (1 - \varepsilon)^{n} = P^{(1)}(E_n).
\]
For a labelled atom $(x, y)$ with $x \neq x_{3}$ and $k \in \{0, 1\}$,
\[
P^{(k)}\{(X, Y) = (x, y) \mid X \neq x_{3}\} = \frac{P(x)\, \mathbf{1}\{y = f_{k}(x)\}}{1 - \varepsilon}.
\]
Since $f_{0}(x) = f_{1}(x)$ for $x \in \{x_{0}, x_{1}, x_{2}\}$, the two conditional laws of one observation coincide. Given $E_n$ the observations are i.i.d.\ with this common law, so $\mathcal{L}_{P^{(0)}}(S_n \mid E_n) = \mathcal{L}_{P^{(1)}}(S_n \mid E_n)$.
\end{proof}

The restriction to the \emph{labelled} observation space is essential: if labels are discarded, the two laws agree everywhere, so the targets remain indistinguishable even after the rare atom is observed.

\begin{proposition}[exact two-point testing risk]\label{prop:two-point}
Let $n \ge 1$, let $\Pi_{0} = (P^{(0)})^{\otimes n}$ and $\Pi_{1} = (P^{(1)})^{\otimes n}$ on the labelled space, and let $R^{*}_n$ be the Bayes risk of the binary test under a uniform prior. Then
\begin{equation}\label{eq:tv-two-point}
\mathrm{TV}(\Pi_{0} ,\ \Pi_{1} ) = 1 - (1 - \varepsilon )^{n},\quad R^{*}_n = \tfrac{1}{2}(1 - \mathrm{TV}(\Pi_{0} ,\Pi_{1} )) = \tfrac{1}{2}(1 - \varepsilon )^{n}.
\end{equation}
\end{proposition}

\begin{proof}
For every event $A \subseteq E_n$ and $k \in \{0, 1\}$, $\Pi_{k}(A) = \Pi_{k}(E_n)\, \mathcal{L}_{P^{(k)}}(A \mid E_n)$, and both factors are the same for $k = 0$ and $k = 1$ by Lemma~\ref{lem:indistinguishability}. Hence $\Pi_{0}$ and $\Pi_{1}$ agree on $E_n$. Every sample in $E_n^{c}$ contains an observation $(x_{3}, y)$, with $y = 0$ under $P^{(0)}$ and $y = 1$ under $P^{(1)}$, so the restrictions of $\Pi_{0}$ and $\Pi_{1}$ to $E_n^{c}$ are concentrated on disjoint sets. Therefore
\[
\mathrm{TV}(\Pi_{0}, \Pi_{1}) = \tfrac{1}{2} \int |d\Pi_{0} - d\Pi_{1}| = \tfrac{1}{2} \bigl(0 + \Pi_{0}(E_n^{c}) + \Pi_{1}(E_n^{c})\bigr) = 1 - (1 - \varepsilon)^{n}.
\]
The Bayes risk of a binary test under a uniform prior is $\tfrac{1}{2}(1 - \mathrm{TV}(\Pi_{0}, \Pi_{1}))$ \cite{lecam1986}, which gives the second identity in Eq.~\eqref{eq:tv-two-point}.
\end{proof}

\begin{proposition}[reduction of family recovery to testing]\label{prop:family-reduction}
For every estimator $\widehat{\mathcal{R}}_n$, including randomised estimators and estimators that know $\varepsilon$ and the class,
\begin{equation}\label{eq:family-lower}
\sup _{P \in \mathcal{P}_{\varepsilon}} P^{n}\{ \widehat{\mathcal{R}}_n \neq \mathcal{R}_P \} \ge \tfrac{1}{2}(1 - \varepsilon )^{n},
\end{equation}
and consequently $n_{\mathrm{family}}(\delta ; \mathcal{P}_{\varepsilon}) \ge \log (1/(2\delta )) / (-\log (1 - \varepsilon )) \ge \log (1/(2\delta )) / (2\varepsilon )$.
\end{proposition}

\begin{proof}
Fix an estimator $\widehat{\mathcal{R}}_n$ and define the test $\psi = \mathbf{1}\{\widehat{\mathcal{R}}_n = \mathcal{R}_{P^{(1)}}\}$, randomised if $\widehat{\mathcal{R}}_n$ is. Since $\mathcal{R}_{P^{(0)}} \neq \mathcal{R}_{P^{(1)}}$ (Table~\ref{tab:populations}),
\[
\{\psi \neq 0\} = \{\widehat{\mathcal{R}}_n = \mathcal{R}_{P^{(1)}}\} \subseteq \{\widehat{\mathcal{R}}_n \neq \mathcal{R}_{P^{(0)}}\}, \qquad \{\psi \neq 1\} = \{\widehat{\mathcal{R}}_n \neq \mathcal{R}_{P^{(1)}}\}.
\]
Hence $P^{(\theta)}\{\psi \neq \theta\} \le P^{(\theta)}\{\widehat{\mathcal{R}}_n \neq \mathcal{R}_{P^{(\theta)}}\}$ for $\theta \in \{0, 1\}$, and by Proposition~\ref{prop:two-point}
\[
\sup_{P \in \mathcal{P}_{\varepsilon}} P^{n}\{\widehat{\mathcal{R}}_n \neq \mathcal{R}_P\} \ge \max_{\theta \in \{0, 1\}} P^{(\theta)}\{\widehat{\mathcal{R}}_n \neq \mathcal{R}_{P^{(\theta)}}\} \ge \tfrac{1}{2} \sum_{\theta = 0}^{1} P^{(\theta)}\{\psi \neq \theta\} \ge R^{*}_n = \tfrac{1}{2}(1 - \varepsilon)^{n}.
\]
For the sample size, $\tfrac{1}{2}(1 - \varepsilon)^{n} \le \delta$ if and only if $n \ge \log (1/(2\delta)) / (-\log (1 - \varepsilon))$, and $-\log (1 - \varepsilon) \le 2\varepsilon$ for $0 < \varepsilon \le 1/2$.
\end{proof}

On the two-point subclass $\{P^{(0)}, P^{(1)}\}$ the bound of Proposition~\ref{prop:family-reduction} is attained. The randomised estimator that reads the label of $x_{3}$ when it appears and otherwise guesses achieves exactly $\tfrac{1}{2}(1-\varepsilon )^{n}$, so on that subclass the family minimax risk equals $\tfrac{1}{2}(1-\varepsilon )^{n}$. On the full class $\mathcal{P}_{\varepsilon}$ only the lower bound is used; the exact minimax value there is not determined.

\subsection{Main theorem}\label{sec:main-theorem}

\begin{theorem}[separation of identification requirements]\label{thm:separation}
Fix $0 < \delta < 1/2$. For the class $\mathcal{P}_{\varepsilon}$ with $0 < \varepsilon \le 1/4$,
\begin{gather}
n_{\mathrm{core}}(\delta ; \mathcal{P}_{\varepsilon}) \le \lceil \log (2/\delta ) / \log (4/3) \rceil \ \text{(independent of}\  \varepsilon ), \label{eq:sep-core}\\
n_{\mathrm{family}}(\delta ; \mathcal{P}_{\varepsilon}) \ge \lceil \log (1/(2\delta )) / (-\log (1 - \varepsilon )) \rceil \ge \log (1/(2\delta )) / (2\varepsilon ), \label{eq:sep-family}
\end{gather}
and consequently
\begin{equation}\label{eq:sep-limit}
\lim _{\varepsilon \downarrow 0} n_{\mathrm{family}}(\delta ; \mathcal{P}_{\varepsilon}) / n_{\mathrm{core}}(\delta ; \mathcal{P}_{\varepsilon}) = \infty .
\end{equation}
In particular, for each fixed $\delta \in (0, 1/2)$ there is no constant $K(\delta )$ such that $n_{\mathrm{family}}(\delta ; \mathcal{P}) \le K(\delta ) \cdot n_{\mathrm{core}}(\delta ; \mathcal{P})$ for every class $\mathcal{P}$ of finite-support realisable populations on which both functionals are non-constant.
\end{theorem}

\begin{proof}[Proof of Theorem~\ref{thm:separation}]

The upper bound is Corollary~\ref{cor:core-uniform}. Since $N(\delta ):= \lceil \log (2/\delta )/\log (4/3)\rceil$ does not involve $\varepsilon$, it is a fixed integer once $\delta$ is fixed. The lower bound is Proposition~\ref{prop:family-reduction}. Both sample-complexity bounds concern minimax risks over the same class $\mathcal{P}_{\varepsilon}$, with the same confidence convention and the same information available to the estimator, so the comparison is between like objects; Proposition~\ref{prop:core-lower} guarantees in addition that the core side is not degenerate.

Dividing,

\begin{equation}\label{eq:sep-ratio}
n_{\mathrm{family}}(\delta ; \mathcal{P}_{\varepsilon}) / n_{\mathrm{core}}(\delta ; \mathcal{P}_{\varepsilon}) \ge \log (1/(2\delta )) / (2 \varepsilon N(\delta )),
\end{equation}

which diverges as $\varepsilon \downarrow 0$ because $\delta$ and hence $N(\delta )$ are fixed. Given $M > 0$, any $\varepsilon < \log (1/(2\delta )) / (2 M N(\delta ))$ forces the ratio above $M$. Integer effects cannot obstruct this: the ceiling only increases the numerator, and the denominator is a fixed integer.
\end{proof}

\emph{On the direction of the bounds.} The upper bound for the core is achieved by the plug-in estimator $\widehat{\mathrm{CORE}}_n$. It is defined for every decision table and uses no knowledge of $\mathcal{P}_{\varepsilon}$, of $\varepsilon$, or of which population generated the data. The lower bound for the family holds for every estimator, including estimators granted full knowledge of the class and of $\varepsilon$. The core upper bound therefore remains valid without class knowledge, whereas the family lower bound remains valid even when that knowledge is supplied.

\emph{The rough-set content.} The statistical step is the standard two-point argument. What is specific to rough sets is Lemma~\ref{lem:discernibility} together with the structural computation of Section~\ref{sec:construction}: a population admits a perturbation of arbitrarily small mass that changes the reduct family while fixing the core. Sections~\ref{sec:atom-addition} and~\ref{sec:witness-cover} isolate that fact, the latter as two population functionals.

\begin{example}[exact plug-in sample sizes]\label{ex:plugin-sizes}
Take $\varepsilon = 10^{-3}$ and $\delta = 0.01$. The plug-in errors on $\mathcal{P}_{\varepsilon}$ are available in closed form. Proposition~\ref{prop:core-exact} gives the core error of each member. By Table~\ref{tab:witnesses}, the family error under $P_{\varepsilon}^{(0)}$ is the probability that one of $x_{0}, x_{2}, x_{3}$ is unobserved, and under $P_{\varepsilon}^{(1)}$ and $P_{\varepsilon}^{(2)}$ it equals the core error. The smallest $n$ at which the worst-case plug-in error is at most $\delta$ is 19 for the core and 4\,603 for the family. The core value coincides with the upper bound of Corollary~\ref{cor:core-uniform}, while the minimax lower bound of Proposition~\ref{prop:core-lower} is 3. The family value exceeds the all-estimator lower bound 3\,911 of Proposition~\ref{prop:family-reduction} by a factor of about $\log (1/\delta )/\log (1/(2\delta )) \approx 1.18$: the plug-in must observe $x_{3}$, whereas an estimator that knows the class may guess its label.
\end{example}

\subsection{Core-preserving atom additions}\label{sec:atom-addition}

The mechanism of Section~\ref{sec:construction}, in which adding one atom keeps the core and changes the family, has an exact local characterisation that does not require recomputing minimal transversals.

\begin{proposition}[atom addition]\label{prop:atom-addition}
Let $P$ be realisable with support $S$,
clutter $\mathcal{C}_P$ and core $\mathrm{CORE}_P$. Let $z$ be a condition vector that does not occur in $S$, with
label $y$, so that the population $Q$ with support $S \cup \{z\}$ and arbitrary positive masses is
again realisable, and let
$N_z = \{ D(z, x): x \in S, f(x) \neq y \}$ be the new edges. Then $\mathrm{CORE}_Q = \mathrm{CORE}_P$ and $\mathcal{R}_Q \neq \mathcal{R}_P$ if and only if both
\begin{gather}
b \in \mathrm{CORE}_P \ \text{for every singleton edge}\ \{b\} \in N_z , \tag{C1}\\
\text{some}\ e \in N_z \ \text{contains no edge of}\ \mathcal{C}_P . \tag{C2}
\end{gather}
\end{proposition}

\begin{proof}
The support of $Q$ is $S \cup \{z\}$, and the conflicting pairs of $Q$ that are not pairs of $P$ are exactly those involving $z$. Hence
\begin{equation*}
\mathcal{H}_Q = \mathcal{H}_P \cup N_z .
\end{equation*}
\emph{Core.} By Lemma~\ref{lem:discernibility}(iii), the singletons of $\mathcal{H}_Q$ are those of $\mathcal{H}_P$ and those of $N_z$, so
\begin{equation*}
\mathrm{CORE}_Q = \mathrm{CORE}_P \cup \{b : \{b\} \in N_z\} .
\end{equation*}
Therefore $\mathrm{CORE}_Q = \mathrm{CORE}_P$ if and only if $\{b : \{b\} \in N_z\} \subseteq \mathrm{CORE}_P$, which is (C1).

\emph{Family.} By Lemma~\ref{lem:discernibility}(ii) and blocker duality \cite{edmonds1970}, $\mathcal{R}_Q \neq \mathcal{R}_P$ if and only if $\mathcal{C}_Q \neq \mathcal{C}_P$. Since every element of $\mathcal{H}_P$ contains an element of $\mathcal{C}_P$,
\begin{equation*}
\mathcal{C}_Q = \min_{\subseteq} (\mathcal{H}_P \cup N_z) = \min_{\subseteq} (\mathcal{C}_P \cup N_z) .
\end{equation*}
\emph{Case (a): every $e \in N_z$ contains an edge of $\mathcal{C}_P$.} Let $e \in N_z$ and $g \in \mathcal{C}_P$ with $g \subseteq e$. No $g' \in \mathcal{C}_P$ satisfies $e \subsetneq g'$, since then $g \subsetneq g'$ would contradict that $\mathcal{C}_P$ is a clutter. So adding $N_z$ removes no edge of $\mathcal{C}_P$ from the minimum, and an element of $N_z$ is minimal only if it equals an edge of $\mathcal{C}_P$; hence $\mathcal{C}_Q = \mathcal{C}_P$.

\emph{Case (b): some $e \in N_z$ contains no edge of $\mathcal{C}_P$.} Choose $e_{0} \in \min_{\subseteq} (\mathcal{C}_P \cup N_z)$ with $e_{0} \subseteq e$. If $e_{0} \in \mathcal{C}_P$, then $e$ would contain an edge of $\mathcal{C}_P$; hence $e_{0} \in \mathcal{C}_Q \setminus \mathcal{C}_P$ and $\mathcal{C}_Q \neq \mathcal{C}_P$.

The two cases are exhaustive, so $\mathcal{R}_Q \neq \mathcal{R}_P$ if and only if (C2) holds. Together with the core part, this proves the proposition.
\end{proof}

The conditions read only the old clutter, the old core, and the new edges. Writing $\mathrm{Agr}(z,x) = C \setminus D(z,x)$, the two conditions translate into support geometry. Condition (C1) says that every oppositely labelled support point at Hamming distance 1 from $z$ differs from $z$ in a core attribute. Condition (C2) says that there are $R \in \mathcal{R}_P$ and an oppositely labelled $x \in S$ with $\pi _R(z) = \pi _R(x)$, that is, a collision on an existing reduct projection. A collision on an old reduct changes the family; (C1) excludes precisely those one-coordinate conflicts that would add a new core attribute.

Proposition~\ref{prop:atom-addition} is close to classical material and serves here as supporting structure. Its ingredients are the discernibility characterisation and blocker duality, and incremental reduct maintenance already updates discernibility information under object insertion \cite{kryszkiewicz1997,wei2018,yang2019}. What is assembled here is the separation of the core condition from the family condition in purely local terms.

Both Proposition~\ref{prop:atom-addition} and its geometric translation were checked by exhaustive enumeration against direct recomputation of the core and the reduct family. The enumeration covered all structures on three binary attributes (34\,752 atom additions) and all structures on four binary attributes with support at most four (828\,800 atom additions). In every case the two criteria agree, and the prediction matches the recomputed structures. The proof does not appeal to these enumerations.

Applied to the construction, it shows how $P_{\varepsilon}^{(0)}$ arises by a single addition. Take $P$ with support $\{x_{0}, x_{1}, x_{2}\}$ and labels $(0,0,1)$, for which $\mathcal{C}_P = \mathcal{R}_P = \{\{a_{1}\}\}$ and $\mathrm{CORE}_P = \{a_{1}\}$. Adding $z = x_{3}$ with label $0$, the only new edge is $\{a_{2},a_{3}\}$, which is not a singleton and does not contain $\{a_{1}\}$; both conditions hold, and the result is $P_{\varepsilon}^{(0)}$.

Proposition~\ref{prop:atom-addition} describes one perturbation at a time. Section~\ref{sec:witness-cover} turns it into statements about populations. There, the smallest masses of atom sets whose deletion changes the family and the core govern both the plug-in error exponents and a non-degenerate class-level separation. For a single atom, those deletions are exactly the reverse of the additions characterised here (Proposition~\ref{prop:single-atom}).

\subsection{Boundary regimes}\label{sec:boundary}

The separation cannot occur everywhere, and the boundary case is worth stating because it recurs empirically in Section~\ref{sec:empirical-populations}.

\begin{definition}[essential attributes]\label{def:essential}
Let the support of $P$ be the full Cartesian product of the attribute value ranges, with labelling $f$. An attribute $a \in C$ is \emph{essential} if two points of the product differing only in the coordinate $a$ carry different labels. $\mathrm{Ess}(f)$ denotes the set of essential attributes.
\end{definition}

\begin{proposition}[full support; folklore in Boolean-function form]\label{prop:full-support}
If the support is the full Cartesian product of the attribute value ranges and the labelling $f$ is non-constant, then
\begin{equation}\label{eq:full-support}
\mathcal{R}_P = \{\mathrm{Ess}(f)\}, \qquad \mathrm{CORE}_P = \mathrm{Ess}(f).
\end{equation}
\end{proposition}

\begin{proof}
Write $\Pi$ for the full Cartesian product, which is the support, and $E = \mathrm{Ess}(f)$.

\emph{Every sufficient $B$ contains $E$.} Let $a \in E \setminus B$. By Definition~\ref{def:essential} there are $x, x' \in \Pi$ with $D(x, x') = \{a\}$ and $f(x) \neq f(x')$. Then $D(x, x') \in \mathcal{H}_P$ and $D(x, x') \cap B = \emptyset$, so $B$ is not sufficient by Lemma~\ref{lem:discernibility}(i).

\emph{$E$ is sufficient.} Let $x, x' \in \Pi$ with $x_{E} = x'_{E}$, and let $b_{1}, \dots , b_{k}$ be the attributes of $C \setminus E$. For $0 \le i \le k$ define $z_{i} \in \Pi$ by $(z_{i})_{b_{j}} = x'_{b_{j}}$ for $j \le i$ and $(z_{i})_{c} = x_{c}$ for all other attributes $c$. Then $z_{0} = x$, $z_{k} = x'$, and $D(z_{i-1}, z_{i}) \subseteq \{b_{i}\}$. Since $b_{i}$ is not essential and every $z_{i}$ lies in $\Pi$, $f(z_{i-1}) = f(z_{i})$, and therefore
\begin{equation*}
f(x) = f(z_{0}) = f(z_{1}) = \dots = f(z_{k}) = f(x') .
\end{equation*}
Hence no two support points that agree on $E$ carry different labels.

By the two claims and the upward closure of sufficiency, $\mathcal{S}_P = \{B \subseteq C : E \subseteq B\}$. Its only inclusion-minimal element is $E$, so $\mathcal{R}_P = \{E\}$ and $\mathrm{CORE}_P = \bigcap \mathcal{R}_P = E$.
\end{proof}

Proposition~\ref{prop:full-support} is the realisable, finite-support counterpart of a familiar probabilistic fact. A distribution satisfying the intersection property, in particular a strictly positive one, has a unique Markov boundary, which consists exactly of the strongly relevant variables \cite{pearl1988,nilsson2007}. Distributions with several Markov boundaries necessarily violate that property \cite{statnikov2013}. Full support plays the role of strict positivity.

On such a population the clutter consists only of singleton edges, so the two plug-in recovery events coincide for every sample, whatever the masses. On any model class consisting entirely of unique-reduct populations, the two targets determine each other and their minimax risks are equal. Proposition~\ref{prop:atom-addition} says the same thing from the other side: a full support offers no hole $z \notin S$ to add.

The converse fails: incomplete support does not imply a non-degenerate separation. Section~\ref{sec:empirical-populations} exhibits a population with incomplete support, multiple reducts, and no meaningful separation at all.

\section{Witness covers and plug-in error exponents}\label{sec:witness-cover}

Theorem~\ref{thm:separation} is a statement about one constructed class. This section asks what governs the two identification problems for an arbitrary population, and answers with two numbers computed from the population.

\subsection{Witness graphs and deletion radii}\label{sec:witness-radii}

Throughout this section $S$ denotes the condition support of $P$; by Assumption~\ref{ass:realisability} each $x \in S$ carries one label, and $P(U) = P(X \in U)$ for $U \subseteq S$.

\begin{definition}[witness graphs and covers]\label{def:witness-graph}
For a nonempty $e \subseteq C$, the \emph{witness graph} $G_e$ has vertex set $S$ and an edge between oppositely labelled support points $x, x'$ with $D(x, x') = e$. A set $K \subseteq S$ meeting every edge of $G_e$ is a \emph{cover} of $G_e$, and $\mu (e)$ denotes the smallest mass $P(K)$ of a cover.
\end{definition}

By Section~\ref{sec:empirical}, $e$ is realised by a sample exactly when some edge of $G_e$ has both endpoints observed.

\begin{definition}[deletion radii and deletion neighbours]\label{def:deletion-radii}
The \emph{family} and \emph{core deletion radii} of $P$ are
\begin{equation}\label{eq:mu-radii}
\mu _{\mathrm{fam}}(P) = \min _{e \in \mathcal{C}_P} \mu (e),\quad \mu _{\mathrm{core}}(P) = \min _{a \in \mathrm{CORE}_P} \mu (\{a\}),
\end{equation}
with $\mu _{\mathrm{fam}}(P) = +\infty$ if $\mathcal{C}_P = \emptyset$ (a single-label support) and $\mu _{\mathrm{core}}(P) = +\infty$ if $\mathrm{CORE}_P = \emptyset$. For $U \subseteq S$ with $P(U) < 1$, the \emph{deletion neighbour} is $P_{\setminus U} = P(\cdot \mid S \setminus U)$.
\end{definition}

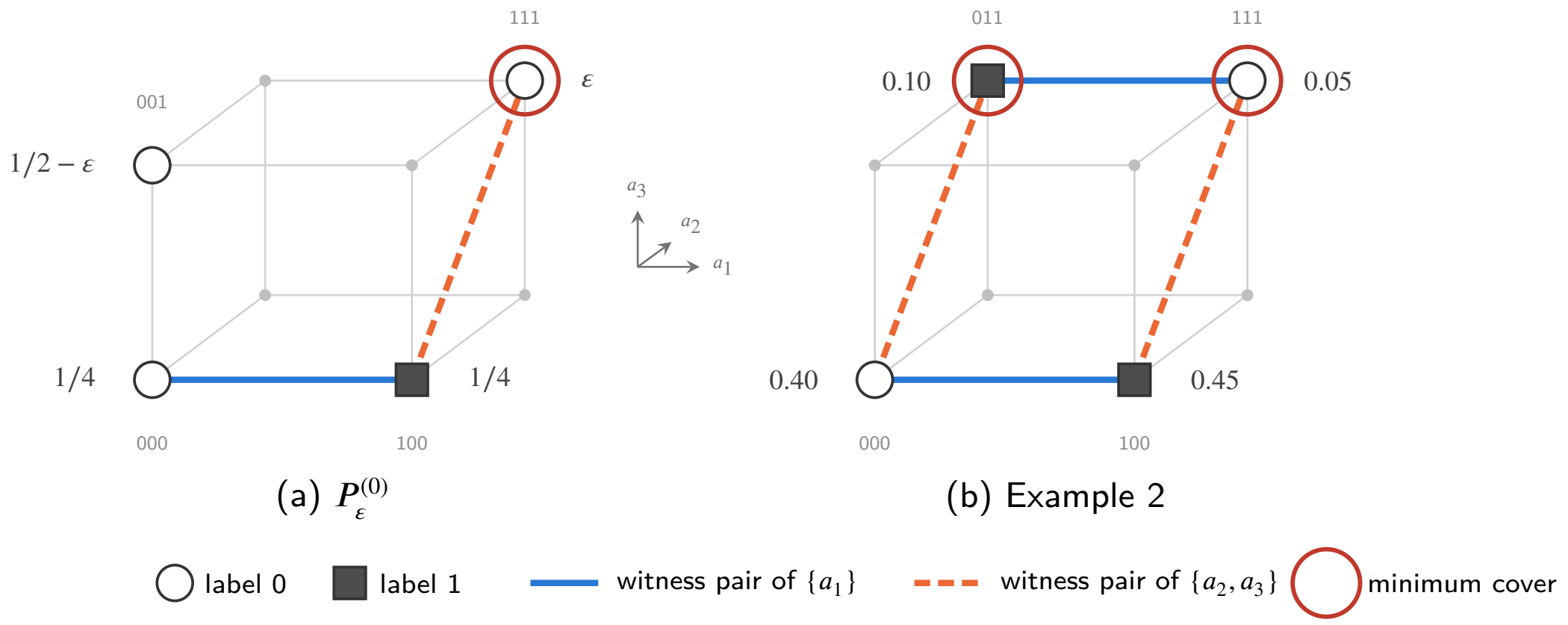
\begin{figure}[pos=!htbp]
\centering
\begin{tikzpicture}[x=1cm,y=1cm,
  lab0/.style={circle,draw=black!80,fill=white,line width=0.7pt,minimum size=9pt,inner sep=0pt},
  lab1/.style={rectangle,draw=black!80,fill=black!70,line width=0.7pt,minimum size=8pt,inner sep=0pt},
  empty/.style={circle,fill=black!25,minimum size=3pt,inner sep=0pt},
  cube/.style={draw=black!18,line width=0.5pt},
  coreedge/.style={draw=wcore,line width=1.6pt},
  famedge/.style={draw=wfam,line width=1.6pt,dash pattern=on 4pt off 2pt},
  cover/.style={circle,draw=wcover,line width=1.1pt,minimum size=17pt,inner sep=0pt},
  mass/.style={font=\scriptsize,text=black!75},
  vlab/.style={font=\tiny,text=black!45}]
\definecolor{wcore}{HTML}{2A78D6}
\definecolor{wfam}{HTML}{EB6834}
\definecolor{wcover}{HTML}{C0392B}
\foreach \pan/\dx in {A/0, B/6.4} {
  \foreach \a/\b/\c in {0/0/0,0/0/1,0/1/0,0/1/1,1/0/0,1/0/1,1/1/0,1/1/1}
    \coordinate (\pan\a\b\c) at ({\dx+2.3*\a+1.0*\b},{1.9*\c+0.75*\b});
  \foreach \u/\v in {000/100,010/110,001/101,011/111,000/010,100/110,001/011,101/111,000/001,100/101,010/011,110/111}
    \draw[cube] (\pan\u) -- (\pan\v);
}
\draw[coreedge] (A000) -- (A100);
\draw[famedge] (A100) -- (A111);
\foreach \v in {010,011,101,110} \node[empty] at (A\v) {};
\node[lab0] at (A000) {}; \node[mass,left=11pt] at (A000) {$1/4$};
\node[lab0] at (A001) {}; \node[mass,left=11pt] at (A001) {$1/2-\varepsilon$};
\node[lab1] at (A100) {}; \node[mass,right=11pt] at (A100) {$1/4$};
\node[lab0] at (A111) {}; \node[mass,right=11pt] at (A111) {$\varepsilon$};
\node[cover] at (A111) {};
\foreach \v/\p in {000/below,001/above,100/below,111/above} \node[vlab,\p=11pt] at (A\v) {\v};
\node[font=\small] at (1.6,-1.05) {(a) $P_{\varepsilon}^{(0)}$};
\draw[coreedge] (B000) -- (B100);
\draw[coreedge] (B011) -- (B111);
\draw[famedge] (B000) -- (B011);
\draw[famedge] (B100) -- (B111);
\foreach \v in {001,010,101,110} \node[empty] at (B\v) {};
\node[lab0] at (B000) {}; \node[mass,left=11pt] at (B000) {$0.40$};
\node[lab1] at (B011) {}; \node[mass,left=11pt] at (B011) {$0.10$};
\node[lab1] at (B100) {}; \node[mass,right=11pt] at (B100) {$0.45$};
\node[lab0] at (B111) {}; \node[mass,right=11pt] at (B111) {$0.05$};
\node[cover] at (B011) {}; \node[cover] at (B111) {};
\foreach \v/\p in {000/below,011/above,100/below,111/above} \node[vlab,\p=11pt] at (B\v) {\v};
\node[font=\small] at (8.0,-1.05) {(b) Example~\ref{ex:covers}};
\begin{scope}[shift={(4.3,1.0)},->,>=stealth,black!55,line width=0.5pt]
  \draw (0,0) -- (0.55,0) node[right,font=\tiny] {$a_1$};
  \draw (0,0) -- (0.30,0.22) node[above right=-1pt,font=\tiny] {$a_2$};
  \draw (0,0) -- (0,0.5) node[above,font=\tiny] {$a_3$};
\end{scope}
\begin{scope}[shift={(0.2,-1.8)},font=\scriptsize]
  \node[lab0] at (0,0) {}; \node[right] at (0.15,0) {label 0};
  \node[lab1] at (1.55,0) {}; \node[right] at (1.7,0) {label 1};
  \draw[coreedge] (3.15,0) -- (3.75,0); \node[right] at (3.8,0) {witness pair of $\{a_1\}$};
  \draw[famedge] (6.55,0) -- (7.15,0); \node[right] at (7.2,0) {witness pair of $\{a_2,a_3\}$};
  \node[cover] at (10.2,0) {}; \node[right] at (10.45,0) {minimum cover};
\end{scope}
\end{tikzpicture}
\caption{Witness graphs on three binary attributes for two populations with clutter $\{\{a_1\}, \{a_2, a_3\}\}$. Grey dots: condition vectors outside the support; numbers: atom masses; ringed atoms: a minimum-mass cover of $G_{\{a_2, a_3\}}$. (a) $\mu_{\mathrm{fam}} = \varepsilon$. (b) $\mu_{\mathrm{fam}} = 0.15$.}
\label{fig:witness}
\end{figure}

The empty-minimum conventions apply only to the deletion radii; the corresponding plug-in error is identically zero. The finite-radius comparisons assume non-empty core. For $e \in \mathcal{H}_P$, $0 < \mu (e) < 1$: a cover of a graph with an edge is nonempty, and $S \setminus \{x\}$ is a cover for every $x$. Singleton edges are minimal, so $\mu _{\mathrm{fam}} \le \mu _{\mathrm{core}}$. A deletion neighbour is again realisable, and $\mathrm{TV}(P, P_{\setminus U}) = P(U)$. Fig.~\ref{fig:witness} shows witness graphs and minimum covers for the construction of Section~\ref{sec:construction} and for the population of Example~\ref{ex:covers} below. Both populations have the same clutter, hence the same reducts and core. Witness pairs of the core edge $\{a_1\}$ differ in one attribute, and those of $\{a_2, a_3\}$ differ in two. In panel (a) each clutter edge has a single witness pair, so the minimum cover of $G_{\{a_2, a_3\}}$ is the rare atom $x_3$ and $\mu_{\mathrm{fam}} = \varepsilon$. In panel (b) each edge has two witness pairs, so a cover needs one atom from each; Example~\ref{ex:covers} computes $\mu_{\mathrm{fam}} = 0.15$, three times the mass of the rarest atom.

The core radius has a familiar interpretation. Group atoms by $z = x_{C \setminus \{a\}}$ and label. For a singleton witness graph, each group is complete multipartite by label, so a minimum cover removes all but a heaviest label class in each group. Consequently,
\begin{gather}
\mu (\{a\}) = \sum _z [P(X_{C \setminus \{a\}} = z) - \max _y P(X_{C \setminus \{a\}} = z,\ Y = y)] = \tau _a, \label{eq:mu-core-tau}\\
\mu _{\mathrm{core}} = \min _{a: \tau _a > 0} \tau _a. \label{eq:mu-core-min}
\end{gather}
Thus $\mu _{\mathrm{core}}$ is the smallest positive leave-one-attribute-out Bayes risk from Section~\ref{sec:population-reducts}. This also gives a polynomial-time computation for any number of labels. The same holds for the family radius. For $e \in \mathcal{C}_P$, group atoms by $z = x_{C \setminus e}$. Two oppositely labelled atoms in one group have a nonempty discernibility set contained in $e$, and minimality of $e$ in $\mathcal{H}_P$ forces it to equal $e$. Hence $G_e$ is again a disjoint union of complete multipartite graphs, and
\begin{equation}\label{eq:mu-edge}
\mu (e) = \sum _z [P(X_{C \setminus e} = z) - \max _y P(X_{C \setminus e} = z,\ Y = y)], \qquad e \in \mathcal{C}_P.
\end{equation}
Every edge of $\mathcal{C}_P$ is the discernibility set of some pair of support points, so $\mathcal{C}_P$ has at most $|S|^2$ edges and $\mu _{\mathrm{fam}}$ is computable in polynomial time for any number of labels.

\subsection{Plug-in error exponents and deletion neighbours}\label{sec:plugin-exponents}

\begin{proposition}[plug-in error exponents]\label{prop:plugin-exponents}
If $\mathcal{C}_P \neq \emptyset$, then for every $n \ge 1$
\begin{equation}\label{eq:plugin-sandwich}
(1 - \mu _{\mathrm{fam}})^{n} \le P^{n}\{ \widehat{\mathcal{R}}_n \neq \mathcal{R}_P \} \le \sum _{e \in \mathcal{C}_P} \sum _K (1 - P(K))^{n},
\end{equation}
the inner sum running over the inclusion-minimal covers $K$ of $G_e$. Write $N_{\mathrm{fam}}$ for the number of summands, so the upper bound is at most $N_{\mathrm{fam}}(1 - \mu _{\mathrm{fam}})^{n}$. In particular $n^{-1} \log P^{n}\{ \widehat{\mathcal{R}}_n \neq \mathcal{R}_P \} \to \log (1 - \mu _{\mathrm{fam}})$. The same holds for $\widehat{\mathrm{CORE}}_n$ with the singleton edges and $\mu _{\mathrm{core}}$; if $\mathrm{CORE}_P = \emptyset$, the plug-in core error is zero.
\end{proposition}

\begin{proof}

Write $M_n = S \setminus O_n$ for the set of unobserved atoms, where $O_n$ is the set of observed atoms. For a \emph{fixed} $K \subseteq S$, every observation avoids $K$ independently with probability $1 - P(K)$, so $P^{n}\{K \subseteq M_n\} = (1 - P(K))^{n}$.

\emph{Event identity.} By Section~\ref{sec:empirical}, $\widehat{\mathcal{R}}_n \neq \mathcal{R}_P$ exactly when some $e \in \mathcal{C}_P$ is unrealised, and $e$ is unrealised exactly when $M_n$ meets every edge of $G_e$, that is, when $M_n$ contains a cover of $G_e$, equivalently an inclusion-minimal cover. Hence

\begin{equation}\label{eq:plugin-event}
\{ \widehat{\mathcal{R}}_n \neq \mathcal{R}_P \} = \bigcup _{e \in \mathcal{C}_P} \bigcup _K \{ K \subseteq M_n \},
\end{equation}

the inner union running over the inclusion-minimal covers $K$ of $G_e$.

\emph{Lower bound.} A cover of mass $\mu _{\mathrm{fam}}$ occurs among the terms, and a union is at least as likely as one of its terms, so the error is at least $(1 - \mu _{\mathrm{fam}})^{n}$.

\emph{Upper bound.} The union bound over the $N_{\mathrm{fam}}$ terms gives $\sum _e \sum _K (1 - P(K))^{n}$, and $P(K) \ge \mu _{\mathrm{fam}}$ for every term, so the sum is at most $N_{\mathrm{fam}} (1 - \mu _{\mathrm{fam}})^{n}$.

\emph{Exponent.} Since $0 < \mu _{\mathrm{fam}} < 1$, taking logarithms in $(1 - \mu _{\mathrm{fam}})^{n} \le P^{n}\{\widehat{\mathcal{R}}_n \neq \mathcal{R}_P\} \le N_{\mathrm{fam}} (1 - \mu _{\mathrm{fam}})^{n}$ and dividing by $n$ gives

\begin{equation}\label{eq:plugin-exponent}
\log (1 - \mu _{\mathrm{fam}}) \le n^{-1} \log P^{n}\{ \widehat{\mathcal{R}}_n \neq \mathcal{R}_P \} \le \log (1 - \mu _{\mathrm{fam}}) + n^{-1} \log N_{\mathrm{fam}},
\end{equation}

and $n^{-1} \log N_{\mathrm{fam}} \to 0$ because $N_{\mathrm{fam}}$ is a finite constant, which proves the limit. The core case is the same argument with the singleton edges of $\mathcal{C}_P$ and $\mu _{\mathrm{core}}$. If $\mathrm{CORE}_P = \emptyset$ then $\widehat{\mathrm{CORE}}_n \subseteq \mathrm{CORE}_P$ forces $\widehat{\mathrm{CORE}}_n = \mathrm{CORE}_P$, so the plug-in core error is zero.
\end{proof}

The probability in Proposition~\ref{prop:plugin-exponents} is elementary: $1 - \mu (e)$ is the largest mass of an independent set of $G_e$, and a one-vertex cover gives the missing-mass probability of a single atom \cite{good1953,mcallester2003}. What is specific to rough sets is that the witness-cover weights determine the exact error exponents of reduct-family and core recovery through $-\log (1 - \mu )$. The constant in the upper bound can be exponential in $|S|$. The radius therefore determines the asymptotic rate, but not by itself an accurate finite-sample threshold. On the finite \texttt{tic-tac-toe} population of Section~\ref{sec:empirical-populations}, $\mu _{\mathrm{fam}} = 2/479$, and at $n = 1\,000$ the lower bound is $0.015$ against an observed family error of $0.19$.

\begin{figure}[pos=!htbp]
\centering
\includegraphics[width=0.72\linewidth]{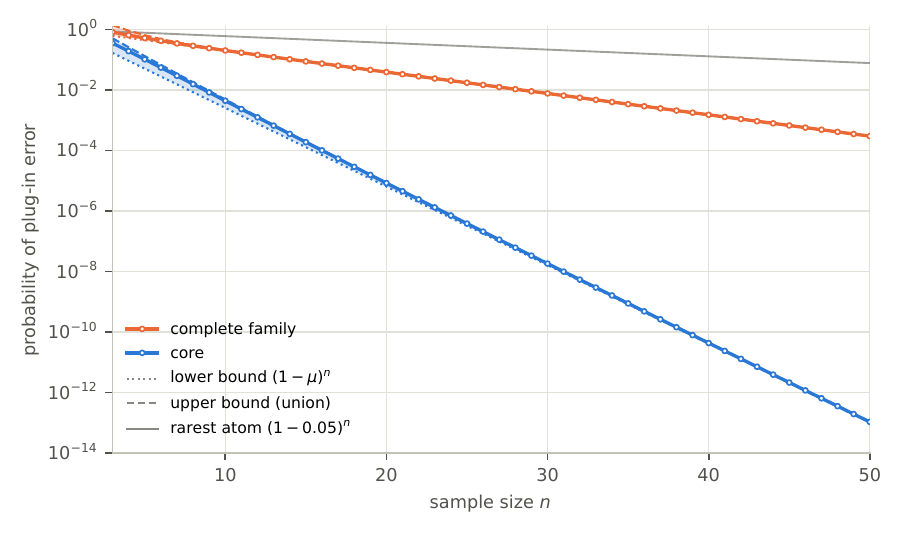}
\caption{Exact plug-in errors for Example~\ref{ex:covers} (markers) with the lower bound $(1-\mu)^n$ (dotted) and the union upper bound (dashed) of Proposition~\ref{prop:plugin-exponents}. Grey line: decay at the rate of the rarest atom.}
\label{fig:example2-bounds}
\end{figure}

\begin{example}[covers that are not single atoms]\label{ex:covers}
Let $S = \{000, 011, 100, 111\}$ on three binary attributes, with labels $0, 1, 1, 0$ (the label is $x_{a_1} + x_{a_2} \bmod 2$) and masses $0.40, 0.10, 0.45, 0.05$. The clutter is $\{\{a_{1}\}, \{a_{2}, a_{3}\}\}$, as for $P_{\varepsilon}^{(0)}$ in Table~\ref{tab:populations}, so $\mathcal{R}_P = \{\{a_{1}, a_{2}\}, \{a_{1}, a_{3}\}\}$ and $\mathrm{CORE}_P = \{a_{1}\}$. Here, however, each clutter edge has two witness pairs (Fig.~\ref{fig:witness}b), and a cover must take one endpoint of each, so $\mu (\{a_{2}, a_{3}\}) = P(\{011, 111\}) = 0.15$ and $\mu (\{a_{1}\}) = P(\{000, 111\}) = 0.45$. The family radius is three times the mass $0.05$ of the rarest atom, because the family changes only when both witness pairs of $\{a_{2}, a_{3}\}$ have lost an atom. The core radius is the leave-one-out Bayes risk of Eq.~\eqref{eq:mu-core-tau}: predicting $Y$ from $(X_{a_2}, X_{a_3})$ errs with probability $0.40$ on the group $00$ and $0.05$ on the group $11$, so $\tau _{a_1} = 0.45$. For the family, Proposition~\ref{prop:plugin-exponents} gives $0.197 \le 0.201 \le 0.205$ at $n = 10$ (lower bound, exact plug-in error, upper bound) and $0.038760 \le 0.038768 \le 0.038775$ at $n = 20$. The smallest $n$ with plug-in error at most $0.01$ is 29 for the family and 9 for the core; Fig.~\ref{fig:example2-bounds} shows the exact errors and the bounds for $3 \le n \le 50$. The family error decays at the rate $-\log(1-\mu_{\mathrm{fam}})$ with $\mu_{\mathrm{fam}} = 0.15$, not at the rate $-\log(1 - 0.05)$ of the rarest atom.
\end{example}

\begin{proposition}[deletion neighbours]\label{prop:deletion}
(a) \emph{(Deletion criterion.)} Let $T$ be a labelled
support, let $U \subseteq T$ with $T \setminus U \neq \emptyset$, and write $\mathcal{R}_T$, $\mathrm{CORE}_T$, $\mathcal{C}_T$ and $G_e[T]$ for the
structures of Sections~\ref{sec:population-reducts} and~\ref{sec:witness-cover} computed on $T$, which depend on the support alone. Then
\begin{gather}
\mathcal{R}_{T \setminus U} \neq \mathcal{R}_T \iff U \ \text{contains a cover of}\ G_e[T] \ \text{for some}\ e \in \mathcal{C}_T , \label{eq:del-fam}\\
\mathrm{CORE}_{T \setminus U} \neq \mathrm{CORE}_T \iff U \ \text{contains a cover of}\ G_{\{a\}}[T] \ \text{for some}\ a \in \mathrm{CORE}_T . \label{eq:del-core}
\end{gather}
Consequently, among proper deletions $U \subseteq S$ with $P(U) < 1$, and with the minimum of an empty set
defined as $+\infty$, $\mu _{\mathrm{fam}}(P)$ is the smallest $P(U)$ with $\mathcal{R}_{P_{\setminus U}} \neq \mathcal{R}_P$, and $\mu _{\mathrm{core}}(P)$ the
smallest over deletions with $\mathrm{CORE}_{P_{\setminus U}} \neq \mathrm{CORE}_P$.
(b) If $K$ is a cover of $G_e$ for some $e \in \mathcal{C}_P$, $P(K) < 1$, and $Q = P_{\setminus K}$, then every estimator $\tilde{\mathcal{R}}_n$, randomised or not, satisfies
\begin{equation}\label{eq:del-lower}
P^{n}\{ \tilde{\mathcal{R}}_n \neq \mathcal{R}_P \} + Q^{n}\{ \tilde{\mathcal{R}}_n \neq \mathcal{R}_Q \} \ge (1 - P(K))^{n}, \qquad \tfrac{1}{2}(1 - P(K))^{n} \le \mathfrak{R}_n(\mathcal{R}; \{P, Q\}) \le (1 - P(K))^{n}.
\end{equation} The same holds for the core when $e = \{a\}$,
$a \in \mathrm{CORE}_P$.
\end{proposition}

\begin{proof}

(a) Deleting atoms only removes conflicting pairs, so $\mathcal{H}_{T \setminus U} \subseteq \mathcal{H}_T$. This inclusion does not transfer to the clutters by itself, because an edge that is not minimal in $\mathcal{H}_T$ can become minimal after smaller edges disappear, so both directions are proved separately.

Suppose first that $U$ contains no cover of $G_e[T]$ for any $e \in \mathcal{C}_T$. Then every $e \in \mathcal{C}_T$ keeps a witness pair, so $\mathcal{C}_T \subseteq \mathcal{H}_{T \setminus U}$; each such $e$ is minimal in $\mathcal{H}_T$, hence minimal in the smaller family $\mathcal{H}_{T \setminus U}$, so $\mathcal{C}_T \subseteq \mathcal{C}_{T \setminus U}$. Conversely every element of $\mathcal{H}_{T \setminus U} \subseteq \mathcal{H}_T$ contains some $e \in \mathcal{C}_T \subseteq \mathcal{H}_{T \setminus U}$, so a minimal element of $\mathcal{H}_{T \setminus U}$ equals that $e$, giving $\mathcal{C}_{T \setminus U} \subseteq \mathcal{C}_T$. The two clutters coincide, and by Lemma~\ref{lem:discernibility}(ii) with blocker duality so do the two families.

Suppose next that $U$ covers $G_e[T]$ for some $e \in \mathcal{C}_T$. Then $e \notin \mathcal{H}_{T \setminus U}$, while $\mathcal{C}_{T \setminus U} \subseteq \mathcal{H}_{T \setminus U}$, so $e$ belongs to $\mathcal{C}_T$ and not to $\mathcal{C}_{T \setminus U}$: the clutters differ, and so do the families.

For the core, $\{a\} \in \mathcal{H}_{T \setminus U}$ implies $\{a\} \in \mathcal{H}_T$, so by Lemma~\ref{lem:discernibility}(iii) a deletion never adds a core attribute; and $a \in \mathrm{CORE}_T$ is lost exactly when no witness pair of $\{a\}$ survives, that is, exactly when $U$ covers $G_{\{a\}}[T]$. Applying this with $T = S$ and reading $P(U)$ as the mass of the deletion (the deletion neighbour $P_{\setminus U}$ has support $S \setminus U$) gives the statements for $P$.

(b) By (a), $\mathcal{R}_Q \neq \mathcal{R}_P$. Write $s$ for a sample sequence. If $s$ avoids $K$ then $Q^{n}(s) = P^{n}(s)/(1 - P(K))^{n} \ge P^{n}(s)$, and otherwise $Q^{n}(s) = 0$, so

\begin{equation}\label{eq:deletion-overlap}
\sum _s \min (P^{n}(s),\ Q^{n}(s)) = \sum _{s \in (S \setminus K)^{n}} P^{n}(s) = (1 - P(K))^{n}.
\end{equation}

Given any estimator $\tilde{\mathcal{R}}_n$, randomised or not, put $\varphi = 1\{\tilde{\mathcal{R}}_n = \mathcal{R}_Q\}$. Under $P$ the event $\{\varphi = 1\}$ forces $\tilde{\mathcal{R}}_n = \mathcal{R}_Q \neq \mathcal{R}_P$, so $\{\varphi = 1\} \subseteq \{\tilde{\mathcal{R}}_n \neq \mathcal{R}_P\}$; under $Q$ the event $\{\varphi = 0\}$ is exactly $\{\tilde{\mathcal{R}}_n \neq \mathcal{R}_Q\}$. The first relation is an inclusion rather than an identity: an estimator returning a third family errs under $P$ while $\varphi$ does not. Hence

\begin{equation}\label{eq:deletion-two-point}
P^{n}\{ \tilde{\mathcal{R}}_n \neq \mathcal{R}_P \} + Q^{n}\{ \tilde{\mathcal{R}}_n \neq \mathcal{R}_Q \} \ge P^{n}\{ \varphi = 1 \} + Q^{n}\{ \varphi = 0 \} \ge \sum _s \min (P^{n}(s),\ Q^{n}(s)) = (1 - P(K))^{n},
\end{equation}

the middle inequality being the two-point argument of Proposition~\ref{prop:family-reduction} applied to the test $\varphi$, and the maximum of the two errors is at least half of $(1 - P(K))^{n}$. For the upper bound, the rule that outputs $\mathcal{R}_P$ exactly when an atom of $K$ has been observed is correct under $Q$, which gives $K$ no mass, and errs under $P$ only when $K$ is unobserved, with probability $(1 - P(K))^{n}$. The core statement is the same argument with $e = \{a\}$, $a \in \mathrm{CORE}_P$.
\end{proof}

Proposition~\ref{prop:deletion} makes the exponents of Proposition~\ref{prop:plugin-exponents} optimal in a precise but limited sense. An estimator that is $\delta$-correct both at $P$ and at its nearest family-changing deletion needs $n \ge \log (1/(2\delta ))/(-\log (1 - \mu _{\mathrm{fam}}))$, while the plug-in is correct at $P$ once $n \ge \lceil \log (N_{\mathrm{fam}}/\delta )/(-\log (1 - \mu _{\mathrm{fam}}))\rceil$. The bounds differ by at most $\log (2N_{\mathrm{fam}})/(-\log (1 - \mu _{\mathrm{fam}})) + 1$, independently of $\delta$ but possibly by a large amount. The statement has two limits. First, it compares requirements at $P$ and says nothing about the accuracy of the plug-in at the neighbour. Second, the neighbourhood consists of deletions only. If a condition vector absent from $S$ may be added with small mass, it can create a new singleton edge, and the core can then be as hard to identify as the family. Deletion neighbours are the appropriate neighbourhood when absent condition vectors are structural zeros of the population, combinations that cannot occur.

\subsection{Class-level separation}\label{sec:class-separation}

Comparing $\mu _{\mathrm{fam}}$ with $\mu _{\mathrm{core}}$ at a single population is a comparison of exponents, not a ratio of minimax sample complexities, which on the one-point class $\{P\}$ both equal one. Turning it into a separation on a common class requires avoiding the degeneracy of Section~\ref{sec:sample-complexity} in a stronger form than non-constancy: the core-changing member of the class must sit at the population's own core radius.

\begin{definition}[core-preserving cover]\label{def:core-preserving}
A cover $K$ of some $G_e$, $e \in \mathcal{C}_P$, is \emph{core-preserving} if every $G_{\{b\}}$, $b \in \mathrm{CORE}_P$, keeps an edge with both endpoints outside $K$.
\end{definition}

\begin{theorem}[class-level separation from witness covers]\label{thm:witness-cover}
Let $\mathrm{CORE}_P \neq \emptyset$ and $\mu _{\mathrm{fam}}(P) < \mu _{\mathrm{core}}(P)$. Let $K$ be a cover of mass $\mu _{\mathrm{fam}}$ of some $G_e$, $e \in \mathcal{C}_P$; let $a \in \mathrm{CORE}_P$ attain $\mu _{\mathrm{core}}$, and let $K_a$ be a cover of $G_{\{a\}}$ of mass $\mu _{\mathrm{core}}$. On the class $\mathcal{K}(P) = \{P, P_{\setminus K}, P_{\setminus K_a}\}$, formed with these choices of $K$ and $K_a$, $K$ is core-preserving, $e$ is not a singleton, both $\mathrm{CORE}$ and $\mathcal{R}$ are non-constant, and for $0 < \delta < 1/2$
\begin{gather}
\log (1/(2\delta )) / (-\log (1 - \mu _{\mathrm{fam}})) \le n_{\mathrm{family}}(\delta ; \mathcal{K}(P)) \le \lceil \log (2/\delta ) / (-\log (1 - \eta ))\rceil , \label{eq:K-family}\\
\log (1/(2\delta )) / (-\log (1 - \mu _{\mathrm{core}})) \le n_{\mathrm{core}}(\delta ; \mathcal{K}(P)) \le \lceil \log (1/\delta ) / (-\log (1 - (\mu _{\mathrm{core}} - \mu _{\mathrm{fam}}))) \rceil . \label{eq:K-core}
\end{gather}
Here $\eta = \min \{\mu _{\mathrm{fam}}, \mu _{\mathrm{core}} - \mu _{\mathrm{fam}}\} > 0$.
\end{theorem}

\begin{proof}

Write $m = \mu _{\mathrm{fam}}(P)$, $m' = \mu _{\mathrm{core}}(P)$, $Q = P_{\setminus K}$, $Q' = P_{\setminus K_a}$, $A = K_a \setminus K$ and $t = P(A)$.

\emph{Step 1 (structures).} If $K$ covered $G_{\{c\}}$ for a core attribute $c$, then $m' \le \mu (\{c\}) \le P(K) = m$, contradicting $m < m'$; the same contradiction rules out a singleton $e$. So $K$ is core-preserving and $e$ is not a singleton. By Proposition~\ref{prop:deletion}(a), $\mathcal{R}_Q \neq \mathcal{R}_P$; a deletion never adds a core attribute, and every core attribute keeps a witness pair outside $K$, so $\mathrm{CORE}_Q = \mathrm{CORE}_P$; and $U = K_a$ covers $G_{\{a\}}[S]$, so $a \notin \mathrm{CORE}_{Q'}$ and $\mathrm{CORE}_{Q'} \neq \mathrm{CORE}_P$. The three families are pairwise distinct: $\mathcal{R}_P \neq \mathcal{R}_Q$ as just shown, and $\mathcal{R}_{Q'}$ differs from both because $\mathrm{CORE}_{Q'} \neq \mathrm{CORE}_P = \mathrm{CORE}_Q$ and the core is the intersection of the family. Both functionals are non-constant on $\mathcal{K}(P)$.

\emph{Step 2 (masses).} The sets $A$ and $K$ are disjoint by construction, and

\begin{equation}\label{eq:K-masses}
t = P(K_a \setminus K) \ge P(K_a) - P(K) = m' - m > 0,\quad Q(A) = P(A)/(1 - m) \ge t,\quad Q'(A) = 0,
\end{equation}

the last equality because $A \subseteq K_a$ and $Q'$ gives $K_a$ no mass.

\emph{Step 3 (lower bounds).} $K$ is a cover of $G_e$ of mass $m < 1$, so Proposition~\ref{prop:deletion}(b) applied to $(P, Q)$ gives $\mathfrak{R}_n(\mathcal{R}; \{P, Q\}) \ge \tfrac{1}{2}(1 - m)^{n}$; enlarging a class cannot decrease the minimax risk, so the same bound holds on $\mathcal{K}(P)$, and $\tfrac{1}{2}(1 - m)^{n} \le \delta$ forces $n \ge \log (1/(2\delta ))/(-\log (1 - m))$. The core bound is Proposition~\ref{prop:deletion}(b) applied to $(P, Q')$ with the singleton edge $\{a\}$ and $P(K_a) = m'$.

\emph{Step 4 (core upper bound).} Let the estimator output $\mathrm{CORE}_{Q'}$ when no atom of $A$ has been observed, and $\mathrm{CORE}_P$ otherwise. Under $Q'$ it is always correct, because $Q'(A) = 0$ leaves $A$ unobserved. Under $P$ it errs exactly when $A$ is unobserved, with probability $(1 - t)^{n}$; under $Q$ exactly when $A$ is unobserved, with probability $(1 - Q(A))^{n} \le (1 - t)^{n}$. The risk over $\mathcal{K}(P)$ is therefore at most $(1 - t)^{n} \le (1 - (m' - m))^{n}$, which is at most $\delta$ once $n \ge \log (1/\delta )/(-\log (1 - (m' - m)))$.

\emph{Step 5 (family upper bound).} Let the estimator output $\mathcal{R}_{Q'}$ when no atom of $A$ has been observed; otherwise $\mathcal{R}_Q$ when no atom of $K$ has been observed, and $\mathcal{R}_P$ otherwise. Under $Q'$ it is always correct. Under $Q$, which gives $K$ no mass, the second branch is the operative one, so it is correct exactly when $A$ is observed and its error is $(1 - Q(A))^{n}$. Under $P$ it is correct exactly when $A$ and $K$ are both observed, so by inclusion--exclusion and $A \cap K = \emptyset$ its error is

\begin{equation}\label{eq:K-family-error}
(1 - t)^{n} + (1 - m)^{n} - (1 - t - m)^{n}.
\end{equation}

With $\eta = \min \{m, m' - m\}$ we have $t \ge m' - m \ge \eta$ and $m \ge \eta$, so all three errors are at most $2(1 - \eta )^{n}$, and $2(1 - \eta )^{n} \le \delta$ once $n \ge \log (2/\delta )/(-\log (1 - \eta ))$.
\end{proof}

Both upper bounds in Theorem~\ref{thm:witness-cover} use estimators that know the class, as minimax risk allows; in Theorem~\ref{thm:separation} the upper bound was attained by the class-agnostic plug-in instead. The core lower bound is what non-constancy alone does not provide: a class whose core-changing member is placed at an arbitrary, heavy cover can satisfy both non-constancy conditions and still have $n_{\mathrm{core}} = 1$. The core lower bound in Eq.~\eqref{eq:K-core} exceeds one when $\mu _{\mathrm{core}} < 1 - 2\delta$.

When $\mu _{\mathrm{core}} \ge 2\mu _{\mathrm{fam}}$, $\eta = \mu _{\mathrm{fam}}$; hence the two family bounds have the same denominator. In particular, for $0 < \delta \le 1/4$,
\begin{equation}\label{eq:K-family-rate}
n_{\mathrm{family}}(\delta ; \mathcal{K}(P)) = \Theta (\log (1/\delta )/\mu _{\mathrm{fam}})
\end{equation}
uniformly in this regime. Indeed $\mu _{\mathrm{fam}} < 1/2$ and $\mu _{\mathrm{fam}} \le -\log (1 - \mu _{\mathrm{fam}}) \le 2\mu _{\mathrm{fam}}$, while the logarithms in the two bounds differ only by universal factors. Without this radius gap, Theorem~\ref{thm:witness-cover} supplies bounds but need not certify a strict sample-complexity gap at a specified $\delta$.

\begin{example}[a separating class from a two-atom cover]\label{ex:witness-class}
For the population of Example~\ref{ex:covers}, $\mu _{\mathrm{fam}} = 0.15 < \mu _{\mathrm{core}} = 0.45$. The minimum cover $K = \{011, 111\}$ of $G_{\{a_2, a_3\}}$ is core-preserving, since it leaves the edge $000$--$100$ of $G_{\{a_1\}}$. Deleting $K$ leaves the population on $\{000, 100\}$, whose only reduct is $\{a_{1}\}$: the family changes and the core does not. Deleting the minimum core cover $K_a = \{000, 111\}$ leaves $\{011, 100\}$, both labelled $1$, so the core becomes empty. On $\mathcal{K}(P)$, Theorem~\ref{thm:witness-cover} gives $25 \le n_{\mathrm{family}} \le 33$ and $7 \le n_{\mathrm{core}} \le 13$ at $\delta = 0.01$, with $\eta = 0.15$; since $\mu _{\mathrm{core}} \ge 2\mu _{\mathrm{fam}}$, the family bounds match in order. Unlike the construction of Section~\ref{sec:construction}, the family-changing deletion removes two atoms, and neither of them alone changes the family.
\end{example}

\begin{corollary}[which supports separate]\label{cor:separating-supports}
(a) A labelled support with non-empty core admits
masses with $\mu _{\mathrm{fam}} < \mu _{\mathrm{core}}$ if and only if some clutter edge has a core-preserving cover.
(b) If masses $P_{\varepsilon}$ on a fixed labelled support with non-empty core satisfy $\mu _{\mathrm{fam}}(P_{\varepsilon}) \le \varepsilon$ and
$\mu _{\mathrm{core}}(P_{\varepsilon}) - \mu _{\mathrm{fam}}(P_{\varepsilon}) \ge \alpha > 0$, then for fixed $\delta \in (0, 1/2)$, with the ratio formed on a single class for each $\varepsilon$,
\begin{equation}\label{eq:sep-supports}
\lim_{\varepsilon \downarrow 0} n_{\mathrm{family}}(\delta ; \mathcal{K}(P_{\varepsilon})) / n_{\mathrm{core}}(\delta ; \mathcal{K}(P_{\varepsilon})) = \infty .
\end{equation}
\end{corollary}

\begin{proof}

(a) If $\mu _{\mathrm{fam}} < \mu _{\mathrm{core}}$, Step 1 of the proof of Theorem~\ref{thm:witness-cover} shows that a cover of mass $\mu _{\mathrm{fam}}$ is core-preserving. Conversely, let $K$ be a core-preserving cover of $G_e$ for some $e \in \mathcal{C}_P$, fix positive masses $p_{0}$ on $S \setminus K$ summing to one, and put $\beta = \min _{x \notin K} p_{0}(x) > 0$. For $0 < \varepsilon \le \varepsilon _{0} = \beta /(2(1 + \beta ))$, give the atoms of $K$ positive masses with total $\varepsilon$ and each $x \notin K$ the mass $(1 - \varepsilon )p_{0}(x)$. Then $\mu _{\mathrm{fam}} \le P(K) = \varepsilon$. Every cover of a core witness graph contains an atom outside $K$, because $K$ is core-preserving, so $\mu _{\mathrm{core}} \ge (1 - \varepsilon )\beta$ and

\begin{equation}\label{eq:cover-gap}
\mu _{\mathrm{core}} - \mu _{\mathrm{fam}} \ge (1 - \varepsilon )\beta - \varepsilon = \beta - \varepsilon (1 + \beta ) \ge \beta /2 = \alpha > 0,
\end{equation}

the last inequality because $\varepsilon \le \beta /(2(1 + \beta ))$. In particular $\mu _{\mathrm{fam}} < \mu _{\mathrm{core}}$, and the hypotheses of (b) hold with this $\alpha$.

(b) Since $\mu _{\mathrm{fam}}(P_{\varepsilon}) \le \varepsilon$ and $\mu _{\mathrm{core}}(P_{\varepsilon}) - \mu _{\mathrm{fam}}(P_{\varepsilon}) \ge \alpha$, the family lower bound and the core upper bound of Theorem~\ref{thm:witness-cover} give

\begin{gather}
n_{\mathrm{family}}(\delta ; \mathcal{K}(P_{\varepsilon})) \ge \log (1/(2\delta )) / (-\log (1 - \varepsilon )), \label{eq:cor-sep-family}\\
n_{\mathrm{core}}(\delta ; \mathcal{K}(P_{\varepsilon})) \le \log (1/\delta ) / (-\log (1 - \alpha )) + 1, \label{eq:cor-sep-core}
\end{gather}

the second using $\lceil u\rceil \le u + 1$. The first bound tends to $+\infty$ as $\varepsilon \downarrow 0$ while the second is a constant determined by $\delta$ and $\alpha$, so the ratio diverges.
\end{proof}

\subsection{Single-atom covers}\label{sec:single-atom}

\begin{proposition}[single atoms]\label{prop:single-atom}
Let $\mathrm{CORE}_P \neq \emptyset$, $z \in S$ and $Q = P_{\setminus \{z\}}$. Then $\{z\}$ is a core-preserving cover of some $G_e$, $e \in \mathcal{C}_P$, if and only if $P$ arises from $Q$ by an atom addition satisfying conditions (C1) and (C2) of Proposition~\ref{prop:atom-addition}, with $Q$ in the role of the original population.
\end{proposition}

\begin{proof}
Since $z \in S$, the population $P$ arises from $Q$ by adding the atom $z$; write $N_z = \{D(z, x) : x \in S \setminus \{z\},\ f(x) \neq f(z)\}$. The two displays in the proof of Proposition~\ref{prop:atom-addition}, with $Q$ in the role of the original population, give
\begin{equation*}
\mathcal{H}_P = \mathcal{H}_Q \cup N_z , \qquad \mathrm{CORE}_P = \mathrm{CORE}_Q \cup \{b : \{b\} \in N_z\} ,
\end{equation*}
so (C1) holds if and only if $\mathrm{CORE}_P = \mathrm{CORE}_Q$. Moreover, for $K = \{z\}$ and $b \in \mathrm{CORE}_P$, the graph $G_{\{b\}}$ keeps an edge avoiding $z$ if and only if $\{b\} \in \mathcal{H}_Q$, that is, $b \in \mathrm{CORE}_Q$. Since $\mathcal{H}_Q \subseteq \mathcal{H}_P$ gives $\mathrm{CORE}_Q \subseteq \mathrm{CORE}_P$, Definition~\ref{def:core-preserving} yields
\begin{equation*}
\{z\} \ \text{is core-preserving} \iff \mathrm{CORE}_Q = \mathrm{CORE}_P \iff \text{(C1)} .
\end{equation*}

\emph{($\Rightarrow$)} Let $\{z\}$ be a core-preserving cover of $G_e$ with $e \in \mathcal{C}_P$. Then (C1) holds by the display above. Every edge of $G_e$ has the endpoint $z$, so every witness pair of $e$ involves $z$; hence $e \in N_z$ and $e \notin \mathcal{H}_Q$. If some $g \in \mathcal{C}_Q$ satisfied $g \subseteq e$, then $g \in \mathcal{H}_Q \subseteq \mathcal{H}_P$ and the minimality of $e$ in $\mathcal{H}_P$ would give $g = e \in \mathcal{H}_Q$, a contradiction. So $e$ contains no edge of $\mathcal{C}_Q$, which is (C2).

\emph{($\Leftarrow$)} Assume (C1) and (C2), and let $e' \in N_z$ contain no edge of $\mathcal{C}_Q$. Since $e' \in \mathcal{H}_P$, there is $e \in \mathcal{C}_P$ with $e \subseteq e'$. If $e \in \mathcal{H}_Q$, then $e$ contains some $g \in \mathcal{C}_Q$, and $g \subseteq e \subseteq e'$ contradicts the choice of $e'$. Hence $e \notin \mathcal{H}_Q$, so every witness pair of $e$ involves $z$ and $\{z\}$ covers $G_e$. By (C1) and the display above, this cover is core-preserving.
\end{proof}

For $P_{\varepsilon}^{(0)}$ with $0 < \varepsilon < 1/4$, the objects of this section are explicit: $\mathcal{C} = \{\{a_{1}\}, \{a_{2},a_{3}\}\}$, $G_{a_{2}a_{3}}$ is the single pair $x_{2}x_{3}$ and $G_{a_{1}}$ the single pair $x_{0}x_{2}$, so $\mu _{\mathrm{fam}} = \varepsilon$ and $\mu _{\mathrm{core}} = 1/4$, with $K = \{x_{3}\}$ and $K_a = \{x_{0}\}$. The class $\mathcal{K}(P_{\varepsilon}^{(0)})$ differs from $\mathcal{P}_{\varepsilon}$, so Theorem~\ref{thm:separation} is an instance of the mechanism of Theorem~\ref{thm:witness-cover} rather than a special case of it.

The results of this section were also checked by exact computation on the finite domains of Table~\ref{tab:checks}. For example, the family estimator of Theorem~\ref{thm:witness-cover} was checked from literal sufficiency and observed-subset laws on 6\,484 labelled supports, 12\,968 mass assignments and 4\,008 risk evaluations across 334 eligible deletion classes, with no violation.

Three statements should be kept distinct: Proposition~\ref{prop:plugin-exponents} gives exact plug-in error exponents; Theorem~\ref{thm:witness-cover} gives minimax bounds on a common deletion class; and Corollary~\ref{cor:separating-supports} characterises supports that admit masses with $\mu _{\mathrm{fam}} < \mu _{\mathrm{core}}$. The last condition does not characterise all separating classes. For example, $P_{\varepsilon}^{(1)}$ has $\mu _{\mathrm{fam}} = \mu _{\mathrm{core}} = 1/4$ but belongs to the separating label-change class $\mathcal{P}_{\varepsilon}$ of Theorem~\ref{thm:separation}.

Which labelled supports admit a core-preserving cover is a combinatorial question about the support; Proposition~\ref{prop:atom-addition} answers it locally for single atoms, and we have no description beyond that. Exhaustive enumeration of small supports shows how restrictive the condition is. Among labelled supports with non-empty core, a core-preserving cover exists for 12.3\% on three binary attributes with at most five atoms and for 41.4\% on four binary attributes with at most four atoms. With three binary attributes, a three-valued decision and at most four atoms, the share is 20.3\%. These are counts over enumerated supports, not a statement about decision tables encountered in practice.

A non-empty core and $\mu _{\mathrm{fam}} < \mu _{\mathrm{core}}$ give the deletion-class bounds of Theorem~\ref{thm:witness-cover}, and arbitrarily large separation follows along the mass sequences of Corollary~\ref{cor:separating-supports}(b), where $\mu _{\mathrm{fam}} \to 0$ while $\mu _{\mathrm{core}} - \mu _{\mathrm{fam}}$ stays bounded below. Neither statement characterises every model class on which separation can occur.

\section{Minimax rates on structural-zero classes}\label{sec:structural-zero}

The classes $\mathcal{P}_{\varepsilon}$ and $\mathcal{K}(P)$ are constructed from a given population and its witness sets, so on a three-point class an upper bound is a test among three known laws. This section replaces them by classes fixed by a modelling assumption. On every such class it determines both minimax rates up to a factor in the risk that depends only on the number $p$ of attributes. The upper bound exceeds the lower bound by at most $4p$ for the core and $4\binom{p}{\lfloor p/2\rfloor}$ for the family.

\subsection{Structural-zero classes and floor radii}\label{sec:sz-classes}

\begin{definition}[structural-zero class]\label{def:sz-class}
Fix a finite set $S$ of condition vectors with a labelling $f$, a set $S_{0} \subseteq S$ of configurations known to occur, and floors $w: S \to (0, 1]$ with $w(S) \le 1$. The \emph{structural-zero class} is
\begin{equation}\label{eq:sz-class}
\mathcal{Q} = \mathcal{Q}(S,\ f,\ S_{0},\ w) = \{ Q : S_{0} \subseteq \operatorname{supp} Q \subseteq S,\ \ \text{labels given by}\  f,\ Q(x) \ge w(x) \ \text{for}\  x \in \operatorname{supp} Q \}.
\end{equation}
A nonempty $T$ with $S_{0} \subseteq T \subseteq S$ is an \emph{admissible support}.
\end{definition}

The class $\mathcal{Q}$ contains every law whose configurations outside $S$ are structural zeros, which contains the configurations in $S_{0}$, and whose present configurations are no rarer than their floors. Every member is realisable, and its structures depend only on its support. Conversely every admissible support $T$ is the support of a member: give each $x \in T$ its floor and one atom of $T$ the remaining mass $1 - w(T) \ge 0$, which is possible because $w(T) \le w(S) \le 1$. The targets $\Psi \in \{\mathcal{R}, \mathrm{CORE}\}$ are therefore functions of the admissible support, and $\Psi$ is \emph{constant on} $\mathcal{Q}$ when it takes the same value on all of them.

For an admissible support $T$ write $\mathcal{H}_T$, $\mathcal{C}_T$, $\mathrm{CORE}_T$ and $G_e[T]$ for the objects of Definitions~\ref{def:reduct}, \ref{def:discernibility} and~\ref{def:witness-graph} computed on $T$.

\begin{definition}[optional covers, changing sets and floor radii]\label{def:floor-radii}
Let $T$ be an admissible support. A cover $V$ of $G_e[T]$ is \emph{optional} if $V \subseteq T \setminus S_{0}$, and a set $U \subseteq T \setminus S_{0}$ is \emph{$\Psi$-changing on} $T$ if $T \setminus U \neq \emptyset$ and $\Psi (T \setminus U) \neq \Psi (T)$. The \emph{floor radii} of $\mathcal{Q}$ are
\begin{gather}
\nu _{\mathrm{fam}} = \min \{ w(V) : S_{0} \subseteq T \subseteq S,\ e \in \mathcal{C}_T,\ V \ \text{an optional cover of}\  G_e[T] \}, \label{eq:nu-fam}\\
\nu _{\mathrm{core}} = \min \{ w(V) : S_{0} \subseteq T \subseteq S,\ a \in \mathrm{CORE}_T,\ V \ \text{an optional cover of}\  G_{\{a\}}[T] \}, \label{eq:nu-core}
\end{gather}
with $\min \emptyset = +\infty$; $\nu _{\Psi}$ denotes the radius of the target $\Psi$.
\end{definition}

Let $U \subseteq T \setminus S_{0}$ with $T \setminus U \neq \emptyset$. By the deletion criterion of Proposition~\ref{prop:deletion}(a), $U$ is $\Psi$-changing on $T$ exactly when it contains an optional cover of $G_e[T]$ for some $e \in \mathcal{C}_T$ if $\Psi = \mathcal{R}$, or of $G_{\{a\}}[T]$ for some $a \in \mathrm{CORE}_T$ if $\Psi = \mathrm{CORE}$. The restriction $T \setminus U \neq \emptyset$ is part of the criterion. It excludes $U = T$, although $T$ covers every witness graph. Each minimum ranges over a finite set, so a finite $\nu _{\Psi}$ is attained by some triple $(T, e, V)$, and then $0 < \nu _{\Psi} < 1$: the floors are positive, and the proof below gives $V \neq T$, so $w(V) = w(T) - w(T \setminus V) < 1$.

\subsection{Minimax rates}\label{sec:sz-rates}

\begin{theorem}[minimax rates on structural-zero classes]\label{thm:structural-zero}
Let $\Psi \in \{\mathcal{R}, \mathrm{CORE}\}$ and let
$\hat{\Psi}_n = \Psi (O_n \cup S_{0})$ be the structure of the set $O_n$ of observed atoms completed by the known configurations. Put $B_p = \binom{p}{\lfloor p/2\rfloor}$, $K_{\mathrm{CORE}} = 2p$ and $K_{\mathcal{R}} = 2B_p$.
(a) $\nu _{\Psi} = +\infty$ if and only if $\Psi$ is constant on $\mathcal{Q}$, and $\nu _{\mathrm{fam}} \le \nu _{\mathrm{core}}$. In both minima the
support may be restricted to $T = S_{0} \cup \{x, x'\}$ for a single oppositely labelled witness pair
$x, x'$ of the edge, and the restricted minimum equals $\nu _{\Psi}$; the two radii are computable in
$O(|S|^{2} (|S_{0}| + 2)^{2} p)$ operations on a unit-cost machine, for any number of decision values.
(b) If $\nu _{\Psi} = +\infty$, then $\hat{\Psi}_n = \Psi (Q)$ for every $Q \in \mathcal{Q}$ and every $n \ge 1$, so $\mathfrak{R}_n(\Psi ; \mathcal{Q}) = 0$ and
$n_{\Psi}(\delta ; \mathcal{Q}) = 1$ for every $\delta \in (0, 1)$. If $\nu _{\Psi} < \infty$, then for every $n \ge 1$
\begin{equation}\label{eq:sz-risk}
\tfrac{1}{2} (1 - \nu _{\Psi})^{n} \le \mathfrak{R}_n(\Psi ; \mathcal{Q}) \le \sup _{Q \in \mathcal{Q}} Q^{n}\{ \hat{\Psi}_n \neq \Psi (Q) \} \le K_{\Psi} (1 - \nu _{\Psi})^{n}.
\end{equation}
Both outer inequalities may be strict. Consequently, for $0 < \delta < 1/2$,
\begin{equation}\label{eq:sz-sample}
\log (1/(2\delta )) / (-\log (1 - \nu _{\Psi})) \le n_{\Psi}(\delta ; \mathcal{Q}) \le \lceil \log (K_{\Psi}/\delta ) / (-\log (1 - \nu _{\Psi})) \rceil .
\end{equation}
\end{theorem}

\begin{proof}

Throughout, $\Psi (T)$ is the structure of the admissible support $T$, and \emph{$\Psi$-changing} is used as defined before the statement.

\emph{Characterisation and reduction.} Let $(T, e, V)$ attain $\nu _{\Psi} < \infty$. Being optional, $V \subseteq T \setminus S_{0}$. If $V = T$ then $T \cap S_{0} = \emptyset$, and for any $x \in T$ the set $T \setminus \{x\}$ is again an optional cover of the same witness graph, because every edge of that graph has two distinct endpoints and so keeps one in $T \setminus \{x\}$; its weight is $w(V) - w(x) < w(V)$ since $w(x) > 0$, contradicting minimality. Hence $V \neq T$, the set $T \setminus V$ is a nonempty admissible support containing $S_{0}$, and $w(V) = w(T) - w(T \setminus V) < 1$. Since $V$ is $\Psi$-changing on $T$, the admissible supports $T$ and $T \setminus V$ satisfy $\Psi (T \setminus V) \neq \Psi (T)$: a finite radius gives a non-constant target.

Conversely, let $T$, $T'$ be admissible supports with $\Psi (T) \neq \Psi (T')$. The admissible support $T \cup T'$ has $\Psi (T \cup T') \neq \Psi (T)$ or $\Psi (T \cup T') \neq \Psi (T')$; in the first case put $U = T' \setminus T$, in the second $U = T \setminus T'$. Either way $U$ is a deletion of $T \cup T'$ that changes $\Psi$ and leaves a nonempty support, and $U$ is optional because $S_{0}$ is contained in both $T$ and $T'$. By the deletion criterion $U$ contains an optional cover, whose weight is finite, so $\nu _{\Psi} < \infty$. For the ordering, every triple feasible for $\nu _{\mathrm{core}}$ is feasible for $\nu _{\mathrm{fam}}$: if $a \in \mathrm{CORE}_T$ then $\{a\} \in \mathcal{C}_T$ by Lemma~\ref{lem:discernibility}(iii), and $G_{\{a\}}[T]$ is the witness graph of that clutter edge; hence $\nu _{\mathrm{fam}} \le \nu _{\mathrm{core}}$.

For the reduction, let $(T, e, V)$ be feasible and pick a witness pair $x, x'$ of $e$ in $T$; put $T' = S_{0} \cup \{x, x'\} \subseteq T$. Then $e = D(x, x') \in \mathcal{H}_{T'} \subseteq \mathcal{H}_T$, and $e$ is minimal in $\mathcal{H}_T$, hence minimal in the smaller family $\mathcal{H}_{T'}$, so $e \in \mathcal{C}_{T'}$; for $\Psi = \mathrm{CORE}$ the edge is a singleton, so its attribute lies in $\mathrm{CORE}_{T'}$. Moreover $V \cap T'$ is an optional cover of $G_e[T']$: every edge of $G_e[T']$ is an edge of $G_e[T]$, hence has an endpoint in $V$, and that endpoint lies in $T'$. Its weight is at most $w(V)$. So the minimum over reduced supports is at most $\nu _{\Psi}$; it is at least $\nu _{\Psi}$ because every reduced triple is itself feasible for the unrestricted minimum. The two minima are therefore equal.

\emph{Computation.} By the reduction it suffices to examine the at most $|S|^{2}/2$ oppositely labelled pairs $x, x'$. For each pair, take the support $T' = S_{0} \cup \{x, x'\}$ of at most $|S_{0}| + 2$ atoms, form $\mathcal{H}_{T'}$ and test whether $e = D(x, x')$ is minimal in it; this costs $O((|S_{0}| + 2)^{2} p)$. Then minimise $w(V)$ over the optional sets $V \subseteq T' \setminus S_{0} \subseteq \{x, x'\}$ covering $G_e[T']$. There are at most four such sets and testing one costs $O((|S_{0}| + 2)^{2})$; if some edge of $G_e[T']$ has both endpoints in $S_{0}$, no optional cover exists and the pair contributes nothing. For $\nu _{\mathrm{core}}$ only pairs with $|D(x, x')| = 1$ are examined. The total is $O(|S|^{2} (|S_{0}| + 2)^{2} p)$ on a unit-cost machine, with no further dependence on the number of decision values. Algorithm~\ref{alg:floor-radii} lists these steps.

\emph{Constant targets and the lower bound.} If $\nu _{\Psi} = +\infty$ then $\Psi$ is constant on $\mathcal{Q}$ by the characterisation above. For $n \ge 1$ the set $O_n \cup S_{0}$ is a nonempty admissible support, so $\hat{\Psi}_n$ equals that constant, which is $\Psi (Q)$ for every member $Q$. The estimator never errs, the minimax risk is zero, and $n_{\Psi}(\delta ; \mathcal{Q}) = 1$.

Let now $\nu _{\Psi} < \infty$ be attained by $(T, e, V)$, and fix $x^{*} \in T \setminus V$, which is nonempty because an attaining cover is proper. Define

\begin{gather}
Q_{1}(x) = w(x) + (1 - w(T)) 1\{x = x^{*}\} \ \text{for}\  x \in T, \label{eq:sz-Q1}\\
Q_{2} = Q_{1}(\cdot \mid T \setminus V). \label{eq:sz-Q2}
\end{gather}

Then $Q_{1}$ is a probability law on $T$ with $Q_{1} \ge w$, so $Q_{1} \in \mathcal{Q}$; and $Q_{1}(T \setminus V) = 1 - w(V) > 0$, so $Q_{2}$ is well defined, satisfies $Q_{2} \ge Q_{1} \ge w$ on $T \setminus V$, and has support $T \setminus V \supseteq S_{0}$, so $Q_{2} \in \mathcal{Q}$. Their targets differ because $V$ is $\Psi$-changing on $T$. A sample $s$ containing an atom of $V$ has $Q_{2}^{n}(s) = 0$; otherwise $Q_{2}^{n}(s) = Q_{1}^{n}(s)/(1 - w(V))^{n} \ge Q_{1}^{n}(s)$. Hence

\begin{equation}\label{eq:sz-overlap}
\sum _s \min (Q_{1}^{n}(s),\ Q_{2}^{n}(s)) = \sum _{s \in (T \setminus V)^{n}} Q_{1}^{n}(s) = (1 - \nu _{\Psi})^{n},
\end{equation}

and the two-point argument of Proposition~\ref{prop:family-reduction}, which covers randomised estimators, bounds the minimax risk over $\{Q_{1}, Q_{2}\}$ below by $\tfrac{1}{2}(1 - \nu _{\Psi})^{n}$. Since $\{Q_{1}, Q_{2}\} \subseteq \mathcal{Q}$, enlarging the class cannot decrease the minimax risk.

\emph{Upper bound and sample complexity.} Fix $Q \in \mathcal{Q}$ with support $T$, and let $E_{\mathcal{R}}(T) = \mathcal{C}_T$ and $E_{\mathrm{CORE}}(T) = \{\{a\} : a \in \mathrm{CORE}_T\}$. The unobserved optional set $U_n = T \setminus (O_n \cup S_{0})$ is random, $T \setminus U_n = O_n \cup S_{0}$ is nonempty, and by the deletion criterion $\hat{\Psi}_n \neq \Psi (Q)$ exactly when $U_n$ contains an optional cover of $G_e[T]$ for some $e \in E_{\Psi}(T)$. For each $e \in E_{\Psi}(T)$ whose witness graph admits an optional cover, fix a witness pair $x_e, x'_e$ of $e$ in $T$ that is not contained in $S_{0}$; one exists, because an edge of $G_e[T]$ with both endpoints in $S_{0}$ cannot be met by an optional set. Put $T_e = S_{0} \cup \{x_e, x'_e\}$. As in the reduction, $e \in \mathcal{C}_{T_e}$, and for the core the attribute of $e$ lies in $\mathrm{CORE}_{T_e}$; moreover, if $V$ is an optional cover of $G_e[T]$, then $V \cap T_e \subseteq \{x_e, x'_e\}$ is an optional cover of $G_e[T_e]$, so the triple $(T_e, e, V \cap T_e)$ is feasible and $w(V \cap T_e) \ge \nu _{\Psi}$. Consequently, if $U_n$ contains an optional cover of $G_e[T]$, it contains an inclusion-minimal optional cover $A \subseteq \{x_e, x'_e\}$ of $G_e[T_e]$. There are at most two such sets $A$, they are fixed in advance, and each has $Q(A) \ge w(A) \ge \nu _{\Psi}$ because $Q \ge w$ on $T$. Writing $\mathcal{A}_e$ for them,

\begin{equation}\label{eq:sz-union}
Q^{n}\{ \hat{\Psi}_n \neq \Psi (Q) \} \le \sum _{e \in E_{\Psi}(T)} \sum _{A \in \mathcal{A}_e} Q^{n}\{ O_n \cap A = \emptyset \} \le 2 |E_{\Psi}(T)| (1 - \nu _{\Psi})^{n}.
\end{equation}

The clutter $\mathcal{C}_T$ is an antichain of subsets of the $p$-element set $C$, so $|\mathcal{C}_T| \le B_p$ by Sperner's theorem, and $|\mathrm{CORE}_T| \le p$; hence $2|E_{\Psi}(T)| \le K_{\Psi}$. Taking the supremum over $Q \in \mathcal{Q}$ gives the third inequality in Eq.~\eqref{eq:sz-risk}, and the minimax risk is at most the risk of this one estimator.

Since $0 < \nu _{\Psi} < 1$, $-\log (1 - \nu _{\Psi}) > 0$. If $\mathfrak{R}_n(\Psi ; \mathcal{Q}) \le \delta$ then the lower bound gives $\tfrac{1}{2}(1 - \nu _{\Psi})^{n} \le \delta$, that is $n \ge \log (1/(2\delta ))/(-\log (1 - \nu _{\Psi}))$, which is the lower bound on $n_{\Psi}(\delta ; \mathcal{Q})$. Conversely $K_{\Psi}(1 - \nu _{\Psi})^{n} \le \delta$ as soon as $n \ge \log (K_{\Psi}/\delta )/(-\log (1 - \nu _{\Psi}))$, and for such an $n$ the upper bound gives $\mathfrak{R}_n(\Psi ; \mathcal{Q}) \le \delta$; the smallest such integer is at most the stated ceiling.
\end{proof}

\subsection{Computing the floor radii}\label{sec:sz-computation}

Algorithm~\ref{alg:floor-radii} computes both floor radii from the class description alone, that is, from $S$, $f$, $S_{0}$ and $w$. It can therefore be run when a study is designed, before any data are collected, and its output gives the sample sizes of Eq.~\eqref{eq:sz-sample}.

\begin{algorithm}[!htbp]
\caption{Floor radii of a structural-zero class}\label{alg:floor-radii}
\begin{algorithmic}[1]
\Require candidate configurations $S$ with labelling $f$, known configurations $S_{0} \subseteq S$, floors $w: S \to (0, 1]$
\Ensure the floor radii $\nu_{\mathrm{fam}}$ and $\nu_{\mathrm{core}}$ of Eqs.~\eqref{eq:nu-fam}--\eqref{eq:nu-core}
\State $\nu_{\mathrm{fam}} \gets +\infty$; \quad $\nu_{\mathrm{core}} \gets +\infty$
\ForAll{pairs $\{x, x'\} \subseteq S$ with $f(x) \neq f(x')$}
  \State $e \gets D(x, x')$; \quad $T' \gets S_{0} \cup \{x, x'\}$ \Comment{reduced support}
  \State $\mathcal{H} \gets \{ D(z, z') : z, z' \in T',\ f(z) \neq f(z') \}$ \Comment{$\mathcal{H}_{T'}$}
  \If{some $h \in \mathcal{H}$ satisfies $h \subsetneq e$} \textbf{continue} \Comment{$e \notin \mathcal{C}_{T'}$}
  \EndIf
  \State $E \gets \{ \{z, z'\} \subseteq T' : f(z) \neq f(z'),\ D(z, z') = e \}$ \Comment{edges of $G_e[T']$}
  \State $c \gets \min \{ w(V) : V \subseteq \{x, x'\} \setminus S_{0},\ V \cap \{z, z'\} \neq \emptyset \ \text{for all}\ \{z, z'\} \in E \}$ \Comment{$\min \emptyset = +\infty$}
  \State $\nu_{\mathrm{fam}} \gets \min (\nu_{\mathrm{fam}}, c)$
  \If{$|e| = 1$} $\nu_{\mathrm{core}} \gets \min (\nu_{\mathrm{core}}, c)$ \Comment{$e = \{a\}$ with $a \in \mathrm{CORE}_{T'}$}
  \EndIf
\EndFor
\State \Return $\nu_{\mathrm{fam}}, \nu_{\mathrm{core}}$
\end{algorithmic}
\end{algorithm}

Each iteration builds at most $(|S_{0}| + 2)^{2}$ discernibility sets of size at most $p$ (lines 4--6), and line 7 tests at most four candidate sets $V$ against at most $(|S_{0}| + 2)^{2}$ edges. Over at most $|S|^{2}/2$ pairs, the running time is $O(|S|^{2} (|S_{0}| + 2)^{2} p)$ and the memory is $O((|S_{0}| + 2)^{2} p)$, as stated in Theorem~\ref{thm:structural-zero}(a). Correctness is the reduction step of the proof: restricting the support to $S_{0} \cup \{x, x'\}$ does not change either minimum.

\subsection{No known configurations}\label{sec:sz-no-known}

When no configuration is known in advance to occur, the radii have a direct geometric reading and the estimator is the class-agnostic plug-in of Section~\ref{sec:empirical}.

\begin{corollary}[no known configurations]\label{cor:no-known-configurations}
Let $S_{0} = \emptyset$ and write $d_H$ for Hamming distance. Then, with $\min \emptyset = +\infty$,
\begin{gather}
\nu _{\mathrm{core}} = \min \{ \min (w(x),\ w(x')) : f(x) \neq f(x'),\ d_H(x,\ x') = 1 \}, \label{eq:nk-nu-core}\\
\nu _{\mathrm{fam}} = \min \{ \min (w(x),\ w(x')) : f(x) \neq f(x') \}. \label{eq:nk-nu-fam}
\end{gather}
If $\nu _{\mathrm{core}} = +\infty$, that is, if no two oppositely labelled configurations lie at Hamming distance one, then $\mathrm{CORE}_Q = \emptyset$ for every $Q \in \mathcal{Q}$, and the plug-in core estimator of Section~\ref{sec:empirical} is exactly correct for every $n \ge 1$. If $\nu _{\mathrm{fam}} = +\infty$, that is, if $f$ is constant on $S$, the same holds for the plug-in family. Each bound below, the case $S_{0} = \emptyset$ of Eq.~\eqref{eq:sz-risk}, applies whenever its own radius is finite:
\begin{equation}\label{eq:nk-plugin}
\sup _{\mathcal{Q}} Q^{n}\{ \widehat{\mathrm{CORE}}_n \neq \mathrm{CORE}_Q \} \le 2p(1 - \nu _{\mathrm{core}})^{n},\quad \sup _{\mathcal{Q}} Q^{n}\{ \widehat{\mathcal{R}}_n \neq \mathcal{R}_Q \} \le 2B_p(1 - \nu _{\mathrm{fam}})^{n},
\end{equation}
with $B_p = \binom{p}{\lfloor p/2\rfloor}$. If $\nu _{\mathrm{core}} < \infty$ and $0 < \delta < 1/2$, then with $L = \log (1/(2\delta ))$, $A_c = \log (2p/\delta )$ and $A_f = \log (2B_p/\delta )$,
\begin{equation}\label{eq:nk-ratio}
\begin{aligned}
&\nu _{\mathrm{core}} L / (2 \nu _{\mathrm{fam}} (A_c + \nu _{\mathrm{core}})) \\
&\quad \le n_{\mathrm{family}}(\delta ; \mathcal{Q}) / n_{\mathrm{core}}(\delta ; \mathcal{Q}) \\
&\quad \le 2 \nu _{\mathrm{core}} (A_f + \nu _{\mathrm{fam}}) / (\nu _{\mathrm{fam}} L).
\end{aligned}
\end{equation}
Both bounds are $\nu _{\mathrm{core}} / \nu _{\mathrm{fam}}$ multiplied by a factor lying in $[L/(2A_c + 1), (2A_f + 1)/L]$, which depends only on $p$ and $\delta$. Consequently, fix $S$, $f$ and $\delta$. Along a sequence of floor profiles with $\nu _{\mathrm{core}} < \infty$, the ratio tends to infinity if and only if $\nu _{\mathrm{fam}}/\nu _{\mathrm{core}} \to 0$. Over a set of such profiles, the ratio is unbounded if and only if $\nu _{\mathrm{fam}}/\nu _{\mathrm{core}}$ has infimum zero on that set. Either way some oppositely labelled pair at Hamming distance at least two must become negligibly rare relative to every oppositely labelled pair at distance one. With a uniform floor $w \equiv \lambda$ and at least one oppositely labelled pair at distance one, $\nu _{\mathrm{fam}} = \nu _{\mathrm{core}} = \lambda$ and the ratio lies in $[1, (2A_f + 1)/L]$.
\end{corollary}

The constant $B_p$ grows exponentially in $p$, but it enters the sample sizes only through $\log B_p$. Table~\ref{tab:p-factors} evaluates the resulting factors at $\delta = 0.01$. For each target, the upper bound in Eq.~\eqref{eq:sz-sample} exceeds the lower bound by at most the factor $\log (K_{\Psi}/\delta ) / \log (1/(2\delta ))$, up to rounding. When $p \le 15$, this factor is below $2.1$ for the core and below $3.6$ for the family, and under a uniform floor the ratio $n_{\mathrm{family}}/n_{\mathrm{core}}$ stays below $7.5$.

\begin{table}[pos=!htbp]
\centering
\small
\caption{Factors of Theorem~\ref{thm:structural-zero} and Corollary~\ref{cor:no-known-configurations} at $\delta = 0.01$ ($L$, $A_c$, $A_f$ as in Corollary~\ref{cor:no-known-configurations}). Columns 3--4: ratio of the upper to the lower sample-size bound; columns 5--6: interval multiplying $\nu _{\mathrm{core}}/\nu _{\mathrm{fam}}$ in Eq.~\eqref{eq:nk-ratio}.}
\label{tab:p-factors}
\begin{tabular}{rrcccc}
\toprule
$p$ & $B_p$ & $\log (2p/\delta)/L$ & $\log (2B_p/\delta)/L$ & $L/(2A_c + 1)$ & $(2A_f + 1)/L$ \\
\midrule
5 & 10 & 1.77 & 1.94 & 0.26 & 4.14 \\
10 & 252 & 1.94 & 2.77 & 0.24 & 5.79 \\
15 & 6\,435 & 2.05 & 3.60 & 0.23 & 7.45 \\
\bottomrule
\end{tabular}
\end{table}

\begin{proof}

\emph{Radii.} By the reduction in the proof of Theorem~\ref{thm:structural-zero} with $S_{0} = \emptyset$, the minima may be taken over supports $T = \{x, x'\}$ with $f(x) \neq f(x')$. On such a $T$ the only discernibility set is $e = D(x, x')$, so $\mathcal{C}_T = \{e\}$, and $G_e[T]$ consists of the single edge $xx'$. Its inclusion-minimal optional covers are the two singletons $\{x\}$ and $\{x'\}$, of weights $w(x)$ and $w(x')$, which gives Eq.~\eqref{eq:nk-nu-fam}. For $\nu _{\mathrm{core}}$, $\mathrm{CORE}_T \neq \emptyset$ requires a singleton clutter edge, that is $d_H(x, x') = 1$, and any such pair gives $\mathrm{CORE}_T = D(x, x')$; the covers are the same, which gives Eq.~\eqref{eq:nk-nu-core}.

\emph{Infinite radii.} If $\nu _{\mathrm{core}} = +\infty$ then $\mathrm{CORE}$ is constant on $\mathcal{Q}$ by Theorem~\ref{thm:structural-zero}(a), and the constant is $\emptyset$ because a one-atom support has empty core. Theorem~\ref{thm:structural-zero}(b) with $S_{0} = \emptyset$, where $\hat{\Psi}_n$ is the plug-in $\widehat{\mathrm{CORE}}_n$ of Section~\ref{sec:empirical}, makes that estimator exactly correct. The family case is identical, and $\nu _{\mathrm{fam}} = +\infty$ says exactly that no two configurations of $S$ carry different labels.

\emph{Constants.} With $S_{0} = \emptyset$ the estimator $\hat{\Psi}_n$ of Theorem~\ref{thm:structural-zero} is the plug-in, and Eq.~\eqref{eq:nk-plugin} is the third inequality of Eq.~\eqref{eq:sz-risk}.

\emph{Ratio.} If $\nu _{\mathrm{core}} < \infty$, it is attained by an oppositely labelled pair $x, x'$ at distance one, and $\min (w(x), w(x')) \le (w(x) + w(x'))/2 \le w(S)/2 \le 1/2$; hence $\nu _{\mathrm{fam}} \le \nu _{\mathrm{core}} \le 1/2$, both radii are finite, and $\nu \le -\log (1 - \nu ) \le 2\nu$ holds for both. The lower bounds of Theorem~\ref{thm:structural-zero}(b) give

\begin{gather}
n_{\mathrm{family}}(\delta ; \mathcal{Q}) \ge L / (-\log (1 - \nu _{\mathrm{fam}})) \ge L / (2\nu _{\mathrm{fam}}), \label{eq:nk-lower-family}\\
n_{\mathrm{core}}(\delta ; \mathcal{Q}) \ge L / (-\log (1 - \nu _{\mathrm{core}})) \ge L / (2\nu _{\mathrm{core}}), \label{eq:nk-lower-core}
\end{gather}

and the plug-in bounds in Eq.~\eqref{eq:nk-plugin}, with $\lceil u\rceil \le u + 1$ and $-\log (1 - \nu ) \ge \nu$, give

\begin{gather}
n_{\mathrm{core}}(\delta ; \mathcal{Q}) \le \lceil A_c / (-\log (1 - \nu _{\mathrm{core}}))\rceil \le A_c/\nu _{\mathrm{core}} + 1 = (A_c + \nu _{\mathrm{core}})/\nu _{\mathrm{core}}, \label{eq:nk-upper-core}\\
n_{\mathrm{family}}(\delta ; \mathcal{Q}) \le (A_f + \nu _{\mathrm{fam}})/\nu _{\mathrm{fam}}. \label{eq:nk-upper-family}
\end{gather}

Dividing the family lower bound by the core upper bound gives the left inequality, and the family upper bound by the core lower bound gives the right one. Written as $\nu _{\mathrm{core}}/\nu _{\mathrm{fam}}$ times a factor, those factors are $L/(2(A_c + \nu _{\mathrm{core}}))$ and $2(A_f + \nu _{\mathrm{fam}})/L$, and $0 < \nu _{\mathrm{fam}} \le \nu _{\mathrm{core}} \le 1/2$ puts them in $[L/(2A_c + 1), (2A_f + 1)/L]$. Since the factors are bounded independently of the profile, both the sequential and the uniform statements follow. Under a uniform floor every oppositely labelled pair carries floor $\lambda$ at both endpoints, so both radii equal $\lambda$; the lower end $1$ of the interval is the coarsening bound $n_{\mathrm{core}} \le n_{\mathrm{family}}$ of Section~\ref{sec:comparison}.
\end{proof}

The construction of Section~\ref{sec:construction} is an instance. On the support of $P_{\varepsilon}^{(0)}$ with floors $1/4$ on $x_{0}, x_{1}, x_{2}$ and $\varepsilon$ on $x_{3}$, the only distance-one conflict is $x_{0}x_{2}$, while $x_{2}x_{3}$ is at distance two; so $\nu _{\mathrm{core}} = 1/4$ and $\nu _{\mathrm{fam}} = \varepsilon$. Theorem~\ref{thm:structural-zero} now applies to a class containing every law with these structural zeros and floors, rather than to three chosen members. Its upper bounds are achieved by the plug-in, and its lower bounds hold for every estimator. The two radii translate directly into sample sizes: at $\delta = 0.01$, Theorem~\ref{thm:structural-zero}(b) gives $14 \le n_{\mathrm{core}} \le 23$ independently of $\varepsilon$, against $194 \le n_{\mathrm{family}} \le 317$ once $\varepsilon = 0.02$. The two intervals differ only by the ratio $-\log(1-1/4)/{-\log(1-\varepsilon)} \approx 14.2$.

\begin{example}[the value of a known configuration]\label{ex:known-config}
Keep the class of the preceding paragraph, with $\varepsilon = 0.02$ and $\delta = 0.01$, and suppose in addition that $x_{3}$ is known to occur, $S_{0} = \{x_{3}\}$. An optional cover must then avoid $x_{3}$, so the witness pair $x_{2}x_{3}$ of $\{a_{2}, a_{3}\}$ can only be covered by $x_{2}$. The family radius rises from $0.02$ to $\nu _{\mathrm{fam}} = 1/4 = \nu _{\mathrm{core}}$, and Theorem~\ref{thm:structural-zero}(b) gives $14 \le n_{\mathrm{family}} \le 23$ in place of $194 \le n_{\mathrm{family}} \le 317$, with the core bounds unchanged. Knowing that the rare configuration occurs, together with its label, removes the separation: the completed estimator $\Psi (O_n \cup S_{0})$ no longer waits for $x_{3}$.
\end{example}

The class $\mathcal{Q}$ contains every admissible sub-support. By Theorem~\ref{thm:structural-zero}(a), both radii, when finite, are therefore attained on a support of the form $S_{0} \cup \{x, x'\}$, and the lower-bound construction lives on such a support. This says nothing about the support of a worst-case law at a given $n$. The floor radii therefore do not record the multiplicity of witnesses measured by the instance radii $\mu _{\mathrm{fam}}$ and $\mu _{\mathrm{core}}$, and Theorem~\ref{thm:structural-zero} describes model classes, while Proposition~\ref{prop:plugin-exponents} and Theorem~\ref{thm:witness-cover} describe a single population. The set $S_{0}$ and the floors are modelling assumptions supplied by the analyst: when $S_{0} = \emptyset$ the plug-in needs neither, but its rate guarantee does. The constants $K_{\Psi}$ depend only on the number of attributes, so the rates are determined uniformly in the size of the support and in $S_{0}$; the family constant can still be exponential in $p$. The role of distance-one pairs mirrors the familiar use of inputs that differ in one coordinate to detect relevant variables. We are not aware of earlier statements of Theorem~\ref{thm:structural-zero} or Corollary~\ref{cor:no-known-configurations}.

Theorem~\ref{thm:structural-zero} and Corollary~\ref{cor:no-known-configurations} were checked by exact computation on the finite domains listed in Appendix~\ref{app:reproducibility}, with code independent of the other checks. They are theoretical results and were not part of the empirical protocol of Section~\ref{sec:empirical-populations}.

\section{Computational checks and an empirical study}\label{sec:experiments}

This section does two things. It checks the exact finite-sample results of Section~\ref{sec:separation} by simulation, and it asks whether the separation appears in finite empirical populations. Neither adds evidence beyond the proofs. Structures are computed exactly by enumeration over attribute subsets, simulations use seeded generators, and Monte Carlo proportions carry pointwise 95\% Wilson intervals \cite{wilson1927}, which measure simulation error only. The candidate datasets, their preprocessing and the computational limits of the empirical study were fixed before any data were downloaded (Appendix~\ref{app:reproducibility}). The witness-cover analysis of Section~\ref{sec:witness-cover} was not part of that study; it explains the outcome afterwards.

\subsection{Validation of the exact results}\label{sec:validation-summary}

\begin{figure}[pos=!htbp]
\centering
\includegraphics[width=0.72\linewidth]{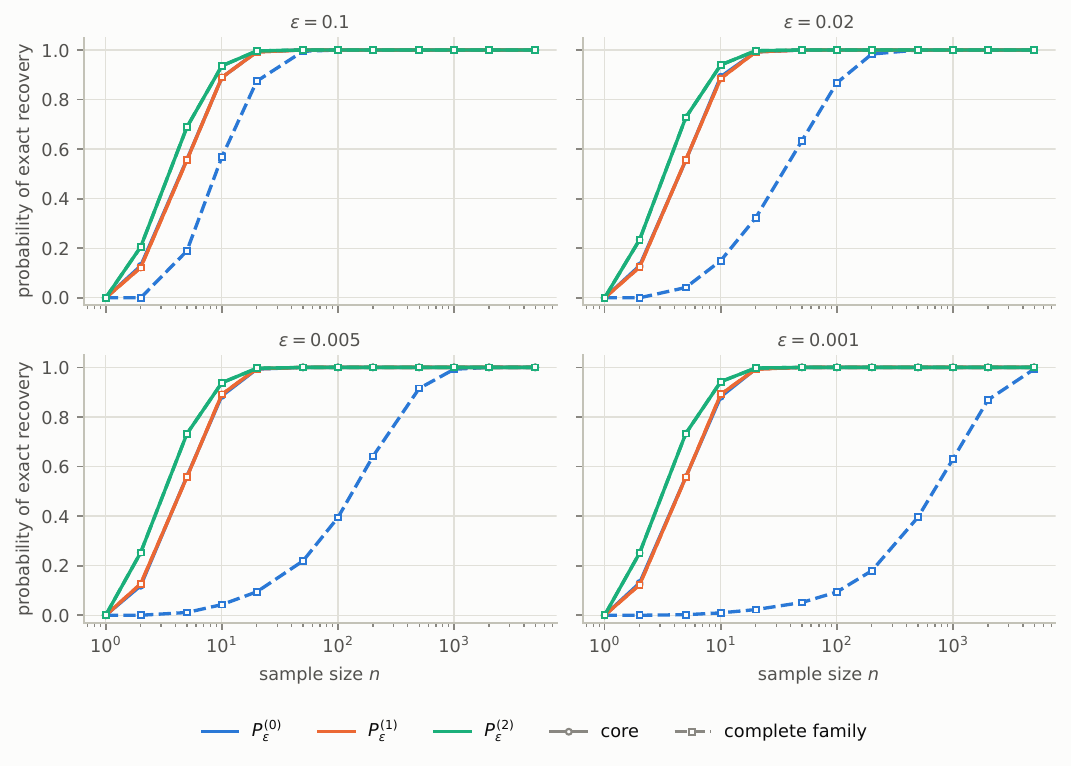}
\caption{Monte Carlo probability of exact plug-in recovery of the core (solid) and the complete family (dashed) on the three-population construction, for four values of $\varepsilon$. Each point: $10^{4}$ or $10^{5}$ repetitions; Wilson half-width at most $0.01$, not drawn.}\label{fig:recovery}
\end{figure}

The simulation checks the recovery probabilities of Proposition~\ref{prop:core-exact} and the plug-in family curves on the three-population construction. At 574 of 576 grid points the Monte Carlo estimates lie within the tolerance of the exact values. The other two points are one estimate, shared by the core and family indicators under $P_{\varepsilon}^{(1)}$, which lies slightly outside the tolerance; a new estimate with $10^{5}$ repetitions agrees with the exact value (Appendix~\ref{sec:exact-validation}). Fig.~\ref{fig:recovery} shows the recovery curves. Only under $P_{\varepsilon}^{(0)}$ does family recovery require the rare atom; under $P_{\varepsilon}^{(1)}$ and $P_{\varepsilon}^{(2)}$ the two recovery events coincide. The waiting times make the separation concrete. The median first sample size at which the empirical core is correct is 5 at every $\varepsilon$ from 0.2 down to 0.001. For the complete family it grows from 7 to 692 and tracks $1/\varepsilon$ (Appendix~\ref{sec:rare-event}).

Fig.~\ref{fig:recovery} compares plug-in probabilities; the separation for every estimator is Theorem~\ref{thm:separation}. Its deterministic bounds at $\delta = 0.01$ are tabulated in Appendix~\ref{sec:sc-separation}: the core bound is 19 for every $\varepsilon$, and the family bound grows as $\log (1/(2\delta ))/\varepsilon$.

\subsection{Finite empirical populations}\label{sec:empirical-populations}

This subsection asks whether the core is identified earlier than the complete family in finite empirical populations beyond the four-atom construction. Each accepted dataset is fixed as a finite empirical population $P_N(x,y) = \#\{i: (X_i,Y_i) = (x,y)\}/N$, and samples are drawn i.i.d.\ with replacement from it, so that $n > N$ is meaningful. $P_N$ is never treated as an unknown real-world law; the ground truth is the exact structure of that finite population.

Of 20 candidate datasets from the UCI Machine Learning Repository \cite{uci}, seven were accepted. Six were excluded because identical condition vectors carry different labels, which violates Assumption~\ref{ass:realisability}. Five exceed the limit $p \le 15$, one is too small, and one exceeds the limit on exact computation. That six of twenty standard categorical decision tables fail Assumption~\ref{ass:realisability} shows that realisability restricts which real tables are in scope.

Recovery probabilities are estimated with $10^{5}$ repetitions at each sample size of the grid
\[
\{0.01, 0.02, 0.05, 0.1, 0.2, 0.5, 1, 2, 5\} \times N \ \cup\ \{1, 2, 5, 10, \dots , 50\,000\},
\]
capped at $5N$. The threshold $n_{\Psi}^{0.95}$ of a target $\Psi$ is the smallest grid size whose Wilson lower bound reaches $0.95$; the point-estimate convention returns the same values on every population. The criterion carries no simultaneous confidence guarantee over the searched grid and no confidence certificate for an unknown population. The descriptive threshold ratio is $\hat{\rho} = n_{\mathrm{family}}^{0.95} / n_{\mathrm{core}}^{0.95}$, and the \emph{coverage} of a population is the fraction of the Cartesian product of the attribute value ranges that occurs in its support. Table~\ref{tab:empirical} collects the accepted populations.

\begin{table}[pos=!htbp]
\centering
\footnotesize\setlength{\tabcolsep}{5pt}
\caption{Accepted finite empirical populations. The ratio for \texttt{tic-tac-toe} is marked degenerate because its core is empty, so the core threshold of $1$ requires no observation.}\label{tab:empirical}
\begin{tabular}{lrrrrrrrrr}
\toprule
dataset & $N$ & $p$ & support & coverage & $|\mathcal{R}_{P_N}|$ & $|\mathrm{CORE}_{P_N}|$ & $n_{\mathrm{core}}^{0.95}$ & $n_{\mathrm{family}}^{0.95}$ & $\hat{\rho}$ \\
\midrule
lenses & 24 & 4 & 24 & 100.0\% & 1 & 4 & 48 & 48 & 1.00 \\
balloons & 20 & 4 & 16 & 100.0\% & 1 & 2 & 20 & 20 & 1.00 \\
balance-scale & 625 & 4 & 625 & 100.0\% & 1 & 4 & 100 & 100 & 1.00 \\
car & 1\,728 & 6 & 1\,728 & 100.0\% & 1 & 6 & 346 & 346 & 1.00 \\
monks-1 & 556 & 6 & 432 & 100.0\% & 1 & 3 & 100 & 100 & 1.00 \\
monks-2 & 601 & 6 & 432 & 100.0\% & 1 & 6 & 100 & 100 & 1.00 \\
tic-tac-toe & 958 & 9 & 958 & 4.9\% & 9 & 0 & 1 & 1\,916 & \emph{degenerate} \\
\bottomrule
\end{tabular}
\end{table}

None of the seven populations shows a non-degenerate separation. Six have a unique reduct and full coverage. For \texttt{monks-1} and \texttt{monks-2}, full coverage follows from the preprocessing, which concatenates the repository's training and test files into the full grid of 432 attribute combinations. On a complete grid the reduct is unique and equals the core (Proposition~\ref{prop:full-support}), every clutter edge is a singleton, and the two recovery events coincide; the measured curves agree at every grid point and $\hat{\rho} = 1.00$. The seventh population, \texttt{tic-tac-toe}, has nine reducts and coverage $4.9\%$, and its thresholds are $1$ for the core against $1\,916$ for the family. Its clutter, however, has no singleton edge, so its core is empty and is recovered from any sample, and the ratio is degenerate. Certificate enumeration was feasible for \texttt{lenses} and \texttt{balloons}, where the core and family certificate sets are identical; for the other five populations the family enumeration exceeded the combination limit.

The negative result does not contradict Theorem~\ref{thm:separation}, which is a possibility statement. The study did not reach the regime in which the question is sharp: incomplete support, multiple reducts and a non-empty core. In the terms of Section~\ref{sec:witness-cover} it could not, because $\mu _{\mathrm{fam}} = \mu _{\mathrm{core}}$ on the six complete grids and \texttt{tic-tac-toe} has no core radius. That the accepted tables fell into these regimes is a limitation of the design, not evidence about observational data.

\section{Discussion}\label{sec:discussion}

This section explains the mechanism behind the results, states their limits and lists the questions they leave open.

\subsection{Support geometry and probability mass}\label{sec:geometry-mass}

How hard each target is to recover depends on the mass of the atoms whose absence can still change it, not on the number of reducts. The witness-cover formulation makes this precise: it identifies the sets of missing observations that can change a target and measures their mass. On structural-zero classes this becomes a statement about support geometry. Every witness pair of a core attribute differs in exactly that attribute. The core is therefore limited by rare conflicts between configurations at Hamming distance one, whereas the family is limited by rare conflicts at any distance.

Under Assumption~\ref{ass:realisability}, the two targets are the family of Markov boundaries and the set of strongly relevant variables (Section~\ref{sec:population-reducts}), so the results can also be read in that language. This reading covers only deterministic, finite-support laws; it does not extend the results beyond them.

Adding an atom only adds conflicting pairs, so its effect can be read from the old clutter and the new edges (Proposition~\ref{prop:atom-addition}). Deleting atoms is not local. An edge disappears only when every one of its witness pairs loses an atom, so the clutter alone is not a sufficient summary; the two populations of Fig.~\ref{fig:witness} share a clutter but not their radii. A sample never loses an atom. Deletion therefore enters only as the lower-bound device of Proposition~\ref{prop:deletion} and as a description of which rare witnesses the sample must see.

On structural-zero classes with no configuration known to occur, Corollary~\ref{cor:no-known-configurations} describes the separation geometrically. Unbounded separation requires oppositely labelled pairs at Hamming distance at least two that become negligibly rare relative to those at distance one. For a single population the corresponding description is open. By Corollary~\ref{cor:separating-supports}, a labelled support with non-empty core admits masses with $\mu _{\mathrm{fam}} < \mu _{\mathrm{core}}$, the sufficient condition of Theorem~\ref{thm:witness-cover}, exactly when some clutter edge has a core-preserving cover. Beyond the single-atom case of Proposition~\ref{prop:atom-addition}, we do not know how to read this condition off the geometry of the support.

\subsection{Limitations}\label{sec:limitations}

All results assume realisability (Assumption~\ref{ass:realisability}), and screening excluded six of twenty candidate datasets for violating it. Nothing here applies to populations with genuine label noise on identical condition vectors. Under such noise the full condition set is no longer sufficient, so the reduct family and the core must first be redefined, for example through positive-region, approximate or variable-precision reducts \cite{ziarko1993,slezak2002}. The one-sidedness of Lemma~\ref{lem:empirical} is also lost: an observed conflicting pair at Hamming distance one may be produced by noise rather than by the population structure, so the empirical core can contain false positives. The waiting-time mechanism behind the lower bounds concerns only unobserved atoms and would persist, but none of the rates is established for such redefined targets.

The results concern exact Pawlak structures. No approximate, probabilistic, fuzzy, neighbourhood or variable-precision extension is claimed, and the recovery events are exact.

The results show that the separation can occur, not how often it occurs. Theorem~\ref{thm:separation} shows that it can be arbitrarily large, and Theorem~\ref{thm:witness-cover} gives a sufficient condition under which a population carries a separating class. Theorem~\ref{thm:structural-zero} gives rates on structural-zero classes once the analyst supplies the structural zeros and floors. None of these results says how typical separating populations are, and no converse to Theorem~\ref{thm:witness-cover} is proved.

The empirical study is narrow. Its datasets are fixed finite populations, sampled with replacement, and no conclusion about real-world populations is drawn from them. Exact computation of the complete family limited the study to $p \le 15$, which excluded five candidates, and certificate enumeration was infeasible for most accepted populations. A large threshold ratio is uninformative when the core is empty, so any study of this ratio should report the core size alongside it.

The instance radii have two limits of scope. First, the exponents $-\log (1 - \mu _{\mathrm{fam}})$ and $-\log (1 - \mu _{\mathrm{core}})$ are sharp for the plug-in estimators and against deletion neighbours only. If new condition vectors may be added with small mass, the core can be as hard to identify as the family. Second, the radii are population quantities. Substituting an empirical population can omit the witness graphs that carry the smallest radius, so the resulting value has no established calibration as a confidence bound for the true radius. Computing the radii is easy, however. Each witness graph of a clutter edge is a disjoint union of complete multipartite graphs, so both radii are computable in polynomial time for any number of decision values (Eq.~\eqref{eq:mu-edge}).

\subsection{Future work}\label{sec:future}

The main structural question is a geometric description of the labelled supports that admit a core-preserving cover, simpler to check than the cover condition itself. Theorem~\ref{thm:witness-cover} gives a matching family rate when $\mu _{\mathrm{core}} \ge 2\mu _{\mathrm{fam}}$, and Theorem~\ref{thm:structural-zero} gives rates on structural-zero classes. Sharper bounds when the two radii are close, the dependence of the family constant on $p$, and classes with unknown floors remain separate questions.

Three extensions need new tools: neighbourhoods that add condition vectors; noisy labels, under which exact sufficient sets may not exist; and sample-based certification or stopping rules. Lemma~\ref{lem:empirical} certifies each observed core attribute individually. Certifying from a sample that no further core attribute or reduct is missing is the main open task for practice.

It is also open whether the realisability screen of Section~\ref{sec:empirical-populations} favours particular support structures, such as complete grids. Complete coverage does not by itself exclude label conflicts: \texttt{monks-3} covers all 432 attribute combinations, yet it was excluded because six condition vectors carry more than one label. Answering the question would require analysing the structure of the excluded datasets.

\section{Conclusion}\label{sec:conclusion}

This paper showed that recovering the core and recovering the whole reduct family from a random sample are statistically different problems. When the analyst can state which configurations may occur and how rare they can be, the difficulty of each target is set by one computable quantity, a weighted cover radius. These radii give the minimax rates up to factors that depend on the number of attributes only (Theorem~\ref{thm:structural-zero}). The core depends only on differently labelled pairs that differ in one attribute, while the family depends on all such pairs. For a known population, the error of the plug-in estimators decays at a rate set by the smallest total mass of atoms whose absence from the sample changes the target. The gap between the two targets has no upper bound. In the simplest example, the number of observations needed for the core does not grow as the rare mass $\varepsilon$ shrinks, while the number needed for the family grows at least like $1/\varepsilon$ (Theorem~\ref{thm:separation}).

In practice, a core attribute found in a sample can be trusted, because an observed conflicting pair proves it. A reduct family found in a sample carries no such proof, since a rare pair that has not yet appeared may still change it. The difference disappears when all combinations of attribute values are present.

All results assume finite Pawlak populations with deterministic labels. They were checked by exact computation and simulation, but seven finite empirical populations gave no example of the separation, so how often it occurs in real data is unknown. A useful next step is a sample-based stopping rule that tells the analyst, from the data alone, when no core attribute or reduct is still missing.

\appendix
\section{Computational evidence and reproducibility}\label{app:reproducibility}

The reproducibility package listed under Data and code availability contains the code, seeds, raw aggregate outputs and processed finite populations behind every number in Section~\ref{sec:experiments} and Appendix~\ref{app:simulation}. These include the Monte Carlo core and family curves, the rare-atom waiting times and the recovery curves of the empirical populations. For the empirical study, it records all 20 candidates of Section~\ref{sec:empirical-populations} with the reasons for exclusion, together with the screening plan, the preprocessing specification and the computational limits. Each processed table has a metadata file with its UCI source URLs, column transformations and value maps; the URLs identify the exact downloaded files. The sample-complexity bounds are reported in Appendix~\ref{sec:sc-separation}.

The package also contains the code and outputs of the exact checks of the theoretical results. Every count quoted in Sections~\ref{sec:separation}--\ref{sec:structural-zero} from an enumeration is reproduced by them, and their scopes are summarised in Table~\ref{tab:checks}.

\begingroup\small\sffamily
\begin{longtable}{@{}p{0.27\linewidth}p{0.40\linewidth}p{0.27\linewidth}@{}}
\caption{Scope and role of the computational checks.}\label{tab:checks}\\
\toprule
Check & Scope & Role \\
\midrule
\endfirsthead
\toprule
Check & Scope & Role \\
\midrule
\endhead
Construction & Three-population structures, exact rational probabilities, testing overlaps and numerical bounds & Independent verification of the construction \\
Atom additions & All labelled supports on three binary attributes (34\,752 atom additions); four binary attributes with at most four atoms (828\,800 atom additions); binary decisions & Proposition~\ref{prop:atom-addition} and its geometric translation \\
Plug-in bounds & All labelled supports on three binary attributes with 2--5 atoms; four binary attributes with 2--4 atoms; three binary attributes and three decision values with 2--4 atoms; specified mass vectors and sample sizes & Plug-in bounds, deletion identities, core-estimator checks and the counts of supports admitting a core-preserving cover (Section~\ref{sec:witness-cover}) \\
Single-atom changes & Single-atom addition/deletion equivalence on its enumerated domains & Proposition~\ref{prop:single-atom} \\
Family upper bound & All labelled supports on three binary attributes with 2--4 atoms and two labels; three binary attributes with 2--3 atoms and three labels; four binary attributes with 2--3 atoms and two labels. Uniform and one seeded rational mass vector per support; $n = 1,\dots ,4$; every minimum changing deletion & Direct sufficiency and observed-subset dynamic programming for the family estimator of Section~\ref{sec:witness-cover}, deletion radii and multiclass core formula \\
Structural-zero classes & All labelled supports on three binary attributes with 2--4 atoms and two labels; sampled supports with three labels (3 attributes, $\le$ 4 atoms), four binary attributes ($\le$ 4 atoms) and three binary attributes ($\le$ 5 atoms); known sets $S_{0}$ and floors drawn per support; Theorem~\ref{thm:witness-cover} estimator on sampled populations & Theorem~\ref{thm:structural-zero} and Corollary~\ref{cor:no-known-configurations}: the constants $K_{\Psi}$, radius by definition, by literal recomputation and by the pair reduction; constancy criterion; exact estimator risks against the upper bounds; lower-bound construction and exact overlaps; Theorem~\ref{thm:witness-cover} error formulas \\
Worked examples & The populations and classes of Examples~\ref{ex:plugin-sizes}--\ref{ex:known-config} & Every number stated in the examples, computed exactly \\
Supplementary sub-claims & Constant targets with infinite radius on the same domains; the range $0 < \nu _{\Psi} < 1$; the support deletion criterion of Proposition~\ref{prop:deletion}(a) over all deletions; the reduction of Theorem~\ref{thm:structural-zero}(a) run as the stated algorithm; the constants of Corollary~\ref{cor:no-known-configurations}; the sandwich of Proposition~\ref{prop:plugin-exponents}; full Cartesian supports for Proposition~\ref{prop:full-support}; the construction of Corollary~\ref{cor:separating-supports}(a) & Targeted checks of specific sub-claims in Propositions~\ref{prop:full-support}--\ref{prop:deletion}, Theorem~\ref{thm:structural-zero} and Corollaries~\ref{cor:separating-supports} and~\ref{cor:no-known-configurations} \\
\bottomrule
\end{longtable}
\endgroup

These finite checks do not enumerate all mass assignments or all estimators. The all-estimator lower bounds follow from the proofs.

\section{Monte Carlo validation of the three-population construction}\label{app:simulation}

This appendix reports the simulation checks summarised in Section~\ref{sec:validation-summary}. Appendix~\ref{sec:exact-validation} compares recovery probabilities with the exact formulas, Appendix~\ref{sec:rare-event} studies waiting times for the rare atom, and Appendix~\ref{sec:sc-separation} evaluates the deterministic sample-complexity bounds.

\subsection{Exact finite-sample validation}\label{sec:exact-validation}

The grid covers $\varepsilon \in \{0.2, 0.1, 0.05, 0.02, 0.01, 0.005, 0.002, 0.001\}$, twelve sample sizes from $n = 1$ to $n = 5\,000$, all three populations and both targets: 576 points. Each point uses $10^{4}$ repetitions. For each population, target and $\varepsilon$, the sample size nearest the $50\%$ crossing of the exact curve uses $10^{5}$, because Monte Carlo noise is largest on this steep part. Both targets are computed from the same simulated samples, so these 48 designated points raise 62 of the 576 grid points to $10^{5}$ repetitions.

The analytic references are themselves verified first. Comparing the closed forms of Proposition~\ref{prop:core-exact} against exhaustive enumeration over observed-atom sets gives 240 identities in floating point and 54 in exact rational arithmetic, with no failure.

Against those references, 574 of 576 Monte Carlo points lie within the tolerance, which is the more permissive of three standard errors and $0.01$. The other two points are one estimate at $\varepsilon = 0.2$, $n = 10$ under $P_{\varepsilon}^{(1)}$, where the core and family witness sets coincide and the two indicators are therefore identical. That estimate, $0.87830$, lies $0.01005$ from the analytic value. A new estimate with $10^{5}$ repetitions and an independent seed gives $0.88814$ against the analytic $0.88835$. Both estimates are kept in the package. The largest standardised deviation over the 214 points with non-degenerate standard error is $3.086$, which is unremarkable for a grid of that size. The remaining 362 points are saturated at 0 or 1, where the standardised deviation carries no information and the absolute error never exceeds $4.4 \cdot 10^{-5}$.

In the simulated core curves, those for $P_{\varepsilon}^{(0)}$ and $P_{\varepsilon}^{(1)}$ coincide and do not depend on $\varepsilon$, as Proposition~\ref{prop:core-exact} requires. The $P_{\varepsilon}^{(2)}$ curve does depend on $\varepsilon$, through $(1/2 + \varepsilon )^{n} - (1/4 + \varepsilon )^{n}$, by about $0.05$ at $n = 2$ between $\varepsilon = 0.1$ and $\varepsilon = 0.001$, while remaining under the $\varepsilon$-free bound $2(3/4)^{n}$ of Corollary~\ref{cor:core-uniform}. The recovery curves are shown in Fig.~\ref{fig:recovery}.

\subsection{The rare-event mechanism}\label{sec:rare-event}

The waiting time for the rare atom, $T_{\mathrm{rare}} = \inf \{n: X_n = x_{3}\}$, is simulated from the four-atom categorical law itself, with $10^{4}$ repetitions per $\varepsilon$ and no censoring at the sample-size cap. Its empirical survival function agrees with $(1 - \varepsilon )^{n}$ within the prescribed tolerance at all 96 grid points ($8$ values of $\varepsilon$ and $12$ sample sizes). A second sampler that draws $T_{\mathrm{rare}}$ directly from a geometric law agrees, but it is reported only as an equivalence check. Verifying $P(T_{\mathrm{rare}} > n) = (1 - \varepsilon )^{n}$ against draws from a law with exactly that survival function would be circular, and the categorical path is what carries the evidence.

The structural part of the experiment runs 2\,000 fully sequential trajectories per $\varepsilon$ under $P_{\varepsilon}^{(0)}$ and records when each structure first becomes correct. By Lemma~\ref{lem:empirical} the recovery time must equal the witness-completion time of Table~\ref{tab:witnesses}, so this part of the experiment is an implementation check rather than evidence. In all 2\,000 trajectories at every $\varepsilon$, the empirical family is correct at the moment the witness set is completed and incorrect one draw earlier, as it must be.

The separation is then visible in the waiting times themselves. The median first sample size at which the empirical core is correct is 5 draws at every $\varepsilon$ from 0.2 down to 0.001. For the complete family the median grows from 7 at $\varepsilon = 0.2$ to 692 at $\varepsilon = 0.001$, tracking $1/\varepsilon$ over more than two orders of magnitude.

Observing $x_{3}$ is necessary for exact family recovery under $P_{\varepsilon}^{(0)}$ but not sufficient. The proportion of trajectories in which the family is already correct at the first appearance of $x_{3}$ rises from $0.394$ at $\varepsilon = 0.2$ to $0.993$ at $\varepsilon = 0.001$. When the rare atom arrives early, the common atoms have often not been seen yet.

\subsection{Sample-complexity separation}\label{sec:sc-separation}

This computation is deterministic. With $n_{\mathrm{core}}^{\mathrm{UB}}(\delta ) = \lceil \log (2/\delta )/\log (4/3)\rceil$ from Corollary~\ref{cor:core-uniform} and $n_{\mathrm{family}}^{\mathrm{LB}}(\delta , \varepsilon ) = \lceil \log (1/(2\delta ))/(-\log (1-\varepsilon ))\rceil$ from Proposition~\ref{prop:family-reduction}, Table~\ref{tab:sc-bounds} gives the values at $\delta = 0.01$, together with $\rho _{\mathrm{LB}} = n_{\mathrm{family}}^{\mathrm{LB}} / n_{\mathrm{core}}^{\mathrm{UB}}$.

\begin{table}[pos=!htbp]
\centering
\small
\caption{Deterministic sample-complexity bounds at $\delta = 0.01$. The last column is the relative error $(\log (1/(2\delta ))/\varepsilon - n_{\mathrm{family}}^{\mathrm{LB}}) / n_{\mathrm{family}}^{\mathrm{LB}}$ of the small-$\varepsilon$ approximation.}\label{tab:sc-bounds}
\begin{tabular}{rrrrrr}
\toprule
$\varepsilon$ & $n_{\mathrm{core}}^{\mathrm{UB}}$ & $n_{\mathrm{family}}^{\mathrm{LB}}$ & $\rho _{\mathrm{LB}}$ & $\log (1/(2\delta ))/\varepsilon$ & rel. error \\
\midrule
$0.1$ & 19 & 38 & 2.0 & 39.1 & $2.9 \cdot 10^{-2}$ \\
$0.01$ & 19 & 390 & 20.5 & 391.2 & $3.1 \cdot 10^{-3}$ \\
$0.001$ & 19 & 3\,911 & 205.8 & 3\,912.0 & $2.6 \cdot 10^{-4}$ \\
$0.0001$ & 19 & 39\,119 & 2\,058.9 & 39\,120.2 & $3.1 \cdot 10^{-5}$ \\
\bottomrule
\end{tabular}
\end{table}

The core column is constant, as Corollary~\ref{cor:core-uniform} requires. A least-squares fit of $\log n_{\mathrm{family}}^{\mathrm{LB}}$ against $\log \varepsilon$ over $10^{-5} \le \varepsilon \le 10^{-2}$, rather than an assumed rate, returns a slope of $-1.0003$, and $\varepsilon \rho _{\mathrm{LB}}$ converges to $0.2059$, which is $\log (1/(2\delta ))/19$, the constant the theorem predicts.

The asymmetry noted in Section~\ref{sec:main-theorem} applies to this table. Here $n_{\mathrm{core}}^{\mathrm{UB}}$ is an upper bound attained by an explicit universal estimator, while $n_{\mathrm{family}}^{\mathrm{LB}}$ is a lower bound over all estimators. $\rho _{\mathrm{LB}}$ is therefore a certified lower bound on the ratio $n_{\mathrm{family}}(\delta ; \mathcal{P}_{\varepsilon}) / n_{\mathrm{core}}(\delta ; \mathcal{P}_{\varepsilon})$, not a ratio of two like-for-like curves.

\section*{CRediT authorship contribution statement}
\textbf{Thanh Dai Tran:} Conceptualization, Methodology, Software, Validation, Formal analysis,
Investigation, Data curation, Writing -- original draft, Writing -- review \& editing, Visualization.

\section*{Declaration of competing interest}
The author declares no known competing financial interests or personal relationships that could have
appeared to influence the work reported in this paper.

\section*{Funding}
This research did not receive any specific grant from funding agencies in the public, commercial, or
not-for-profit sectors.

\section*{Data and code availability}
The finite empirical populations of Section~\ref{sec:empirical-populations} are derived from publicly available
datasets of the UCI Machine Learning Repository \cite{uci}; the exact source file of each dataset is recorded in the reproducibility package. The
reproducibility package described in Appendix~\ref{app:reproducibility}, containing code, seeds, raw aggregate outputs,
processed populations, checksums, the verification scripts for the experimental results of
Section~\ref{sec:experiments} and Appendix~\ref{app:simulation} and the exact checks summarised in Table~\ref{tab:checks}, is available at
\url{https://github.com/ttdaiuneti/rough-learnability-reproducibility}.

\bibliographystyle{cas-model2-num-unsrt}
\bibliography{references}

\end{document}